# Nonlinear Cruise Controllers with Bidirectional Sensing for a String of Vehicles

Iasson Karafyllis[*], Dionysios Theodosis[**], and Markos Papageorgiou[**,***]

[*]Dept. of Mathematics, National Technical University of Athens,
Zografou Campus, 15780, Athens, Greece,
emails: iasonkar@central.ntua.gr , iasonkaraf@gmail.com

[**]Dynamic Systems and Simulation Laboratory,
Technical University of Crete, Chania, 73100, Greece
(emails: dtheodosis@dssl.tuc.gr , markos@dssl.tuc.gr)

[***]Faculty of Maritime and Transportation,
Ningbo University, Ningbo, China

## Abstract

We introduce a nonlinear cruise controller that is fully decentralized (by vehicle) and uses spacing and speed measurements from the preceding and following vehicles to decide on the appropriate control action (acceleration) for each vehicle. The proposed cruise controller is studied on both a ring-road and an open road and guarantees that there are no collisions between vehicles, while their speeds are always positive and never exceed the road speed limits. For both cases of the open road and the ring-road, we rigorously prove that the set of equilibrium points is globally asymptotically stable and provide *KL* estimates that guarantee uniform convergence to the said set. Moreover, we show that for the ring-road, and under certain conditions, there is a single equilibrium point which is exponentially attractive.

## 1. Introduction

Over the past few decades, driver-assistance systems have been incorporated in vehicles that provide a safer and more comfortable driving experience. The conventional Cruise Control (of vehicle speed) and the Adaptive Cruise Control (ACC) (of vehicle speed and of distance to the front vehicle) are such systems, having the potential to increase safety, reduce traffic accidents, and improve traffic flow on highways (see [10], [13], [17], [21], and references therein). Recent technological developments have also introduced more advanced systems for autonomous driving, like the Cooperative Adaptive Cruise Control (CACC) which utilizes vehicle-to-vehicle communications and may increase safety and efficiency (see for instance [3], [27], [38]). While most proposed or deployed cruise controllers rely on measurements of the downstream vehicle (distance, speed), an alternative approach that provides various advantages, is the bidirectional cruise control system, first proposed in [36], where each vehicle monitors both the preceding and the following vehicles, see also [1], [4], [9], [11], [16], [19], [23], [24], [31], [35], [38].

Of great academic interest are the so-called ring-roads, where the first vehicle tracks the speed and position of the last vehicle (see for instance [6], [12], [22], [25], [26], [37]). Several traffic phenomena on a ring-road may give important insights about traffic jams that also appear on open roads (see [2], [30], [33]). Also, several recent studies have shown the ability of connected and automated vehicles on ring-roads to dissipate traffic waves, see [7], [8], [15].

From a mathematical point of view, the design of a cruise controller is a particularly challenging problem, even if one assumes vehicle motion in only one spatial dimension (1-D motion or lane-based motion) because: (i) the control system is defined on an open set (the state space); (ii) the controller must be completely distributed (per vehicle), using measurements from neighbouring vehicles only; and (iii) the set of equilibria to be stabilized may be an unbounded set. The last requirement has important consequences: for systems with non-compact invariant sets, we cannot use standard concepts like the concept of a size function (see [28]) or standard results for attractors (see [5]). Such problems are rarely studied in control theory. The geometry of the state space and the set of equilibria for the cruise controller design problem is unique in control theory, as a neighbourhood of the set of equilibria can contain the whole state space. We are not aware of other control problems that have all the above characteristics.

In this paper we present a novel cruise controller that has the following main features: (a) bidirectional sensing: the control action (vehicle acceleration) uses the distance and relative speeds from its preceding and following vehicles; (b) the cruise controller imposes low measurement demands, using data only when the distance between vehicles is within a specified range, while there is no interaction between vehicles when the distance is sufficiently large; (c) the vehicles do not collide with each other; (d) the speeds of all vehicles are always positive and remain below a road speed limit; (e) all vehicle speeds converge to the (common) desired speed; and (f) all the above features hold globally, i.e., for all physically relevant initial conditions. We rigorously prove that those features hold for both the case of an open road and the case of a ring-road. Due to the various state constraints discussed above, the overall system of automated vehicles on both the ring-road and the open road is defined in open sets. Thus, we need advanced tools to study the properties of such systems. Specifically, we employ size functions and size-like functions (see [28]) combined with Lyapunov functions to show the following qualitative properties, for the ring-road and open road, respectively.

<u>Ring-road</u>: (i) We show global asymptotic stability with respect to two measures of the set of equilibria and establish a *KL* estimate for the solutions of the closed-loop system that guarantees uniform convergence rate to the set of equilibria (Theorem 2). (ii) Under certain conditions on the length of the ring-road and the measurement demands of the vehicles, we show that there exists a single equilibrium point that is locally exponentially stable and globally exponentially attractive (Theorem 3). The latter is a consequence of a novel result for dynamical systems defined in open sets (systems with state constrains) which states that local exponential stability and uniform global asymptotic stability imply global exponential convergence (Theorem 1).

<u>Open road</u>: (i) We show global asymptotic stability with respect to two measures of the set of equilibria and establish a *KL* estimate for the solutions of the closed-loop system that guarantees uniform convergence rate to the set of equilibria (Theorem 4). (ii) We also show that the inter-vehicle distance is always bounded and depends on the initial spacing and speed of the vehicles (Proposition 2). This fact allows us also to show that the closed-loop system is Lagrange stable (Corollary 1). Finally, we show that, under certain conditions, the vehicles' speeds converge exponentially fast to the desired speed (Proposition 3).

Hardly any cruise controllers in the literature can guarantee all the above features. The cruise controller in [17] adopts the Follow-the-Leader architecture and leads the vehicles to a platoon formation. Consequently, the set of equilibria is a compact set (a singleton). Moreover, the state space in [17] is different from the one in the present work and involves relations between inter-vehicle distances and vehicle speeds. In [16], a bidirectional cruise controller presented with analogous properties as above. However, there are several important differences with the controller presented in this paper. Here, forward completeness (i.e., the fact that the solution remains in the state space for all times) is guaranteed directly by using an appropriate Lyapunov function, while in [17] forward completeness follows from an aggressive adjustment of the controller gains. Moreover, here the desired speed is not constant, but depends on the distance from the preceding and following vehicles. For instance, the desired speed may decrease when the distance to the preceding and following vehicles is small.

It should be highlighted that the proposed cruise controller is based on the 2-D (lane-free) cruise controller proposed in [18] when restricted to lane-based motion. However, there are some essential differences. Specifically, we are using here measurements only from the preceding and the following vehicle while the lane-free controller in [18] requires measurements from all vehicles within an interaction distance. Another difference is that the proposed lane-based controller contains a non-zero "viscosity" term which uses relative speed measurements from the preceding and the following vehicles. More importantly, we show that in lane-based driving, the proposed cruise controller has several qualitative properties that are not available in the lane-free case.

The structure of the paper is as follows. Section 2 is devoted to the problem formulation. Section 3 addresses certain properties of dynamical systems defined in open sets. Section 4 presents the cruise controllers with bidirectional architecture and the main results for the case of the ring-road. Section 5 presents the properties of the proposed cruise controller for the case of an open road. All proofs of the main results are provided in Section 6. Section 7 presents numerical examples to demonstrate the properties of the cruise controller and comparisons with the controller presented in [16]. Finally, some concluding remarks are given in Section 8.

**Notation.** Throughout this paper, we adopt the following notation. $\mathbb{R}_+ := [0,+\infty)$ is the set of non-negative real numbers. By $|x|$ we denote both the Euclidean norm of a vector $x \in \mathbb{R}^n$ and the absolute value of a scalar $x \in \mathbb{R}$. By $1_n$ we denote the 1-vector of size $n$, i.e., $1_n = (1,1,...,1) \in \mathbb{R}^n$. By $C^0(A,\Omega)$, we denote the class of continuous functions on $A \subseteq \mathbb{R}^n$, which take values in $\Omega \subseteq \mathbb{R}^m$. By $C^k(A;\Omega)$, where $k \geq 1$ is an integer, we denote the class of functions on $A \subseteq \mathbb{R}^n$ with continuous derivatives of order $k$, which take values in $\Omega \subseteq \mathbb{R}^m$. When $\Omega = \mathbb{R}$ the we write $C^0(A)$ or $C^k(A)$. For a set $S \subseteq \mathbb{R}^n$, $\overline{S}$ denotes the closure of $S$. By $dist(x,A)$ we denote the Euclidean distance of the point $x \in \mathbb{R}^n$ from the set $A \subset \mathbb{R}^n$, i.e., $dist(x,A) = \inf\{|x-y|: y \in A\}$. By $K$ we denote the class of increasing $C^0$ functions $a: \mathbb{R}_+ \to \mathbb{R}_+$ with $a(0) = 0$. By $K_\infty$ we denote the class of increasing $C^0$ functions $a: \mathbb{R}_+ \to \mathbb{R}_+$ with $a(0) = 0$ and $\lim_{s \to +\infty} a(s) = +\infty$. By $KL$ we denote the set of all continuous functions $\sigma: \mathbb{R}_+ \times \mathbb{R}_+ \to \mathbb{R}_+$ with the following properties: (i) for each $t \geq 0$ the mapping $\sigma(\cdot,t)$ is of class $K$; (ii) for each $s \geq 0$, the mapping $\sigma(s,\cdot)$ is non-increasing with $\lim_{t \to +\infty} \sigma(s,t) = 0$.

## 2. Problem Formulation

The longitudinal movement of $n \geq 2$ vehicles on a ring-road is described by the following set of ODEs:

$$\begin{aligned} \dot{s}_i &= v_{i-1} - v_i \quad , \quad i = 1,...,n \\ \dot{v}_i &= F_i \quad , \quad i = 1,...,n \end{aligned} \quad (2.1)$$

with $v_0 = v_n$ and $\sum_{i=1}^{n} s_i = R$, where $R > 0$ is the length of the ring-road, where $s_i$ is the back-to-back distance between vehicle $i$ and its preceding vehicle $i-1$, $v_i$ is the speed of vehicle $i$, and $F_i$ is the acceleration of vehicle $i$. To achieve collision avoidance, we require that for all times $t \geq 0$, the inter-vehicle distances $s_i(t)$, $i = 1,...,n$, should be greater than a certain safety distance $L > 0$. In addition, we want to ensure that vehicles never move backwards, i.e., $v_i(t) > 0$, respect the speed limit of the road $v_{max} > 0$, i.e., $v_i(t) < v_{max}$ for all $t \geq 0$, and eventually tend to the same desired speed $v^* \in (0, v_{max})$.

System (2.1) is equivalently described by the system

$$\begin{aligned} \dot{s}_i &= v_{i-1} - v_i \quad , \quad i = 2,...,n \\ \dot{v}_i &= F_i \quad , \quad i = 1,...,n \end{aligned} \tag{2.2}$$

with state space the open set

$$\Omega = \left\{ \begin{array}{l} (s_2,...,s_n, v_1,...,v_n) \in \mathbb{R}^{2n-1}: \\ \min_{i=2,...,n}(s_i) > L, L + \sum_{i=2}^{n} s_i < R \\ \max_{i=1,...,n}(v_i) < v_{max}, \min_{i=1,...,n}(v_i) > 0 \end{array} \right\} \tag{2.3}$$

where $L \in (0, n^{-1}R)$. Notice that the state space of model (2.2) is an open and bounded set.

Analogously, we can also consider $n \geq 2$ vehicles on an open road whose longitudinal movement is described by (2.2) with the state-space being an open and unbounded set given by

$$\Omega_S = \left\{ \begin{array}{l} (s_2,...,s_n, v_1,...,v_n) \in \mathbb{R}^{2n-1}: \\ \min_{i=2,...,n}(s_i) > L, \\ \max_{i=1,...,n}(v_i) < v_{max}, \min_{i=1,...,n}(v_i) > 0 \end{array} \right\}. \tag{2.4}$$

The state-spaces $\Omega$ and $\Omega_S$ described by (2.3) and (2.4), for the ring-road and the open road, respectively, act as a basis for the design of cruise controllers and for expressing the main objectives of the paper:

- Each vehicle uses only the distance and relative speed from its preceding and following vehicles to select the proper control action (vehicle acceleration);
- The cruise controllers are fully decentralized and do not require measurements from neighboring vehicles whose distance is greater that a certain threshold;
- the vehicles do not collide with each other;
- the speeds of all vehicles are always non-negative and remain below an a priori given speed limit $v_{max}$;
- all vehicle speeds converge (asymptotically and/or exponentially) to a given longitudinal speed set-point $v^* \in (0, v_{max})$;
- all inter-vehicle distances converge to an equilibrium position;
- all the above features are valid globally, i.e., hold for all physically relevant initial conditions.

Before proceeding with the systematic design of the cruise controllers for the ring-road and the open road, we will study some properties of dynamical systems defined in open sets (systems with state constraints).

# 3. Dynamical Systems Defined in Open Sets

Consider the dynamical system

$$\dot{x} = f(x)$$
$$x \in D \subseteq \mathbb{R}^n \tag{3.1}$$

where $D \subseteq \mathbb{R}^n$ is an open set and $f: D \to \mathbb{R}^n$ is a locally Lipschitz mapping. We first recall the definition of a size function for (3.1), see also [28].

**Definition: (Size function for (3.1)):** *A continuous function $\varphi: D \to \mathbb{R}_+$ that satisfies the following property:*

*"For every $r > \inf(\varphi(D))$ the set $\{x \in D : \varphi(x) \leq r\}$ is a compact set"*

*is called a size function for (3.1).*

It is clear that a continuous function $\varphi: D \to \mathbb{R}_+$ is a size function for (3.1) if and only if it satisfies the following property:

"For every $r > \inf(\varphi(D))$ there exists a compact set $K \subseteq D$

such that $\{x \in D : \varphi(x) \leq r\} \subseteq K$"

We next assume that $0 \in D$ is an equilibrium point for (3.1), i.e., $f(0) = 0$. We have the following result.

**Theorem 1 (local exponential stability and uniform global asymptotic stability implies global exponential convergence):** *Suppose that $0 \in D$ is locally exponentially stable for (3.1), i.e., suppose that there exist constants $\rho, \omega > 0$ and $M \geq 1$ such that for every $x_0 \in D$ with $|x_0| \leq \rho$, the solution $x(t)$ of the initial-value problem (3.1) with $x(0) = x_0$ satisfies the estimate*

$$|x(t)| \leq M \exp(-\omega t)|x_0|, \text{ for all } t \geq 0 \tag{3.2}$$

*Moreover, suppose that $0 \in D$ is uniformly globally asymptotically stable, i.e., there exists $\sigma \in KL$ and a size function $\varphi: D \to \mathbb{R}_+$ for (3.1) with $\varphi(0) = 0$ such that for every $x_0 \in D$ the solution $x(t)$ of the initial-value problem (3.1) with $x(0) = x_0$ satisfies the estimate*

$$|x(t)| \leq \sigma(\varphi(x_0), t), \text{ for all } t \geq 0 \tag{3.3}$$

*Then there exists a size function $\kappa: D \to \mathbb{R}_+$ for (3.1) with $\kappa(0) = 0$ such that for every $x_0 \in D$ the solution $x(t)$ of the initial-value problem (3.1) with $x(0) = x_0$ satisfies the estimate*

$$|x(t)| \leq \exp(-\omega t)\kappa(x_0), \text{ for all } t \geq 0 \tag{3.4}$$

**Proof:** By virtue of Proposition 7 in [29] there exist functions $a_1, a_2 \in K_\infty$ such that

$$\sigma(s,t) \leq a_1(\exp(-t)a_2(s)), \text{ for all } t, s \geq 0 \tag{3.5}$$

We also define:

$$\beta(s) := s + \sigma(s,0), \text{ for all } s \geq 0 \tag{3.6}$$

$$T(x) = \ln\left(\max\left(1, \frac{a_2(\varphi(x))}{a_1^{-1}(\rho)}\right)\right), \text{ for all } x \in D \tag{3.7}$$

Notice that $\beta \in K_\infty$.

We next prove the following claim: For every $x_0 \in D$, the solution $x(t)$ of the initial-value problem (3.1) with $x(0) = x_0$ satisfies the estimate

$$|x(t)| \le M \exp(-\omega(t - T(x_0)))\beta(\varphi(x_0) + |x_0|). \tag{3.8}$$

In order to prove (3.8), we distinguish the following cases.

Case 1: $|x_0| \le \rho$. In this case, (3.2) holds. Using (3.6) and (3.7) (which imply $\exp(\omega T(x_0))\beta(\varphi(x_0) + |x_0|) \ge |x_0|$), it follows that (3.8) holds.

Case 2: $|x_0| > \rho$. In this case, (3.3) holds. Using (3.3), (3.5) and (3.7)

$$|x(T(x_0))| \le \sigma(\varphi(x_0), T(x_0)) \le a_1\left(\exp(-T(x_0))a_2(\varphi(x_0))\right) \le \rho. \tag{3.9}$$

Using the semigroup property, it follows from (3.2) that the following estimate holds:

$$|x(t)| \le M \exp(-\omega(t - T(x_0)))|x(T(x_0))|, \text{ for all } t \ge T(x_0). \tag{3.10}$$

Exploiting (3.10), (3.3) and (3.6), we get for $t \ge T(x_0)$:

$$\begin{aligned}
|x(t)| &\le M \exp(-\omega(t - T(x_0)))\sigma(\varphi(x_0), T(x_0)) \\
&\le M \exp(-\omega(t - T(x_0)))\sigma(\varphi(x_0), 0) \\
&\le M \exp(-\omega(t - T(x_0)))\beta(\varphi(x_0)) \\
&\le M \exp(-\omega(t - T(x_0)))\beta(\varphi(x_0) + |x_0|)
\end{aligned} \tag{3.11}$$

On the other hand, exploiting estimate (3.3) and the fact that $M \ge 1$, we get for $t \in [0, T(x_0)]$:

$$\begin{aligned}
|x(t)| &\le \sigma(\varphi(x_0), t) \le \sigma(\varphi(x_0), 0) \\
&\le \beta(\varphi(x_0)) \le M \exp(-\omega(t - T(x_0)))\beta(\varphi(x_0)) \\
&\le M \exp(-\omega(t - T(x_0)))\beta(\varphi(x_0) + |x_0|)
\end{aligned} \tag{3.12}$$

Consequently, (3.8) holds in every case. Estimate (3.4) is a direct consequence of (3.8) with

$$\kappa(x) = M \exp(\omega T(x))\beta(\varphi(x) + |x|) = M\left(\max\left(1, \frac{a_2(\varphi(x))}{a_1^{-1}(\rho)}\right)\right)^\omega \beta(\varphi(x) + |x|),$$

Notice that the fact that $\varphi : D \to \mathbb{R}_+$ is a size function for (3.1) with $\varphi(0) = 0$ and the fact that $\beta \in K_\infty$ guarantees that $\kappa : D \to \mathbb{R}_+$ is a size function for (3.1) with $\kappa(0) = 0$. The proof is complete. ◁

We next give a technical lemma that generalizes Proposition 2.2 on page 107 in [14].

**Lemma 1:** *Let $D \subset O \subset \mathbb{R}^n$ be given non-empty sets and assume that $O \subset \mathbb{R}^n$ is open. Let also $Q, W : O \to \mathbb{R}_+$ be continuous functions with $Q(x) = W(x) = 0$ for all $x \in D$ and $Q(x) > 0$ for all $x \in O \setminus D$. Moreover, suppose that $Q(O) = \mathbb{R}_+$ and that for every $r \geq 0$ the set $\{x \in O : Q(x) \leq r\}$ is compact. Then there exists a function $\zeta \in K_\infty$ such that $W(x) \leq \zeta(Q(x))$, for all $x \in O$. Moreover, if $W(x) > 0$ for all $x \in O \setminus D$ then there exists a positive definite and globally Lipschitz function $\rho : \mathbb{R}_+ \to \mathbb{R}_+$ such that*

$$W(x) \geq \rho(Q(x)), \text{ for all } x \in O. \tag{3.13}$$

**Proof:** Define

$$\bar{\zeta}(s) := \max\{W(x) : x \in O, Q(x) \leq s\} \tag{3.14}$$

Definition (3.14) guarantees that $\bar{\zeta} : \mathbb{R}_+ \to \mathbb{R}_+$ is non-decreasing with $\bar{\zeta}(0) = 0$. Moreover, definition (3.14) implies that $W(x) \leq \bar{\zeta}(Q(x))$, for all $x \in O$.

We next claim that $\lim_{s \to 0^+}(\bar{\zeta}(s)) = 0$. The proof is made by means of a contradiction argument. Suppose the contrary, i.e., that there exists $\varepsilon > 0$ and a sequence $\{b_k > 0 : k = 1, 2, ...\}$ with $\lim_{k \to +\infty}(b_k) = 0$ such that $\bar{\zeta}(b_k) \geq \varepsilon$ for $k = 1, 2, ...$. It follows from (3.14) that there exists a sequence $\{x_k \in O : k = 1, 2, ...\}$ with $Q(x_k) \leq b_k$ and $W(x_k) \geq \varepsilon$ for $k = 1, 2, ...$. Notice that the sequence $\{x_k \in O : k = 1, 2, ...\}$ is bounded, since $\{x_k \in O : k = 1, 2, ...\} \subseteq \{x \in O : Q(x) \leq \sup_{k \geq 1}(b_k)\}$. Consequently, there exists a convergent subsequence of the sequence $\{x_k \in O : k = 1, 2, ...\}$. Therefore, without loss of generality we may assume that the sequence $\{x_k \in O : k = 1, 2, ...\}$ is convergent, i.e., there exists $x^* \in \{x \in O : Q(x) \leq \sup_{k \geq 1}(b_k)\}$ with $\lim_{k \to +\infty}(x_k) = x^*$. Since $0 \leq Q(x_k) \leq b_k$, $\lim_{k \to +\infty}(b_k) = 0$ and $W(x_k) \geq \varepsilon$ for $k = 1, 2, ...$, by continuity of $Q, W : O \to \mathbb{R}_+$ we get $Q(x^*) = 0$ and $W(x^*) \geq \varepsilon$; a contradiction.

Lemma 2.4 on page 65 in [14] implies the existence of $\zeta \in K_\infty$ such that $\bar{\zeta}(s) \leq \zeta(s)$, for all $s \geq 0$. Thus, $W(x) \leq \zeta(Q(x))$, for all $x \in O$.

Next, we assume that if $W(x) > 0$ for all $x \in O \setminus D$. Notice that our assumptions guarantee that the sets $\{x \in O : 1 \leq Q(x) \leq s\}$ for $s > 1$ and $\{x \in O : s \leq Q(x) \leq 1\}$ for $0 \leq s \leq 1$ are non-empty and compact. Define

$$\bar{\rho}(s) := \min\{W(x) : x \in O, 1 \leq Q(x) \leq s\}, \text{ for } s > 1 \tag{3.15}$$

and

$$\bar{\rho}(s) := \min\{W(x) : x \in O, s \leq Q(x) \leq 1\}, \text{ for } 0 \leq s \leq 1. \tag{3.16}$$

Definitions (3.15) and (3.16) imply that $\bar{\rho}$ is positive definite. Moreover, definitions (3.15) and (3.16) also imply that $\bar{\rho}$ is non-decreasing on $[0,1]$, non-increasing on $[1, +\infty)$ and satisfies $W(x) \geq \bar{\rho}(Q(x))$ for all $x \in O$. Finally, define the function

$$\rho(s) = \inf\{\bar{\rho}(y) + |y - s| : y \geq 0\}, \; s \geq 0. \tag{3.17}$$

By standard inf-convolution arguments, it follows that $\rho$ defined by means of (3.17), is positive definite, globally Lipschitz and satisfies $0 \leq \rho(s) \leq \bar{\rho}(s)$ for all $s \geq 0$. Thus, we get that $W(x) \geq \rho(Q(x))$ for all $x \in O$. This concludes the proof. ◁

## 4. The Cruise Controller for the Ring-Road

To design cruise controllers with bidirectional sensing and collision avoidance, we employ artificial potential functions, see [34]. We require the cruise controller to be fully decentralized, relying on measurements of distance and speed of both the preceding and following vehicles, only when their distance is less than the interaction distance $\lambda > L$. Let $V \in C^3((L,\lambda); \mathbb{R}_+)$ be a function that satisfies

$$\lim_{x \to L^+} (V(x)) = +\infty, \tag{4.1}$$

$$V(x) = 0, \quad x \geq \lambda. \tag{4.2}$$

$$\begin{aligned} V'(x) &< 0, \quad \text{for } x \in (L, \lambda) \\ V''(x) &> 0, \quad \text{for } x \in (L, \lambda) \end{aligned}. \tag{4.3}$$

Properties (4.1) and (4.3) imply that the potential $V(s)$ exerts a repulsion force when vehicles are close to each other, while property (4.2) implies that there is no interaction when vehicles are distant. In what follows, we use the convention $s_{n+1} = s_1 = R - \sum_{i=2}^{n} s_i$, $v_{n+1} = v_1$. Define

$$s = (s_2, \ldots, s_n) \in \mathbb{R}^{n-1}, \quad v = (v_1, \ldots, v_n) \in \mathbb{R}^n \tag{4.4}$$

and the Lyapunov function

$$H(s,v) := \frac{v_{\max}^2}{2} \sum_{i=1}^{n} \frac{(v_i - f_i(s))^2}{v_i(v_{\max} - v_i)} + \sum_{i=1}^{n} V(s_i) \tag{4.5}$$

where

$$f_i(s) = v^* - b(V'(s_{i+1}) - V'(s_i)), i = 1, \ldots, n \tag{4.6}$$

$v^* \in (0, v_{\max})$ is the desired speed of the vehicles, and $b: \mathbb{R} \to (v^* - v_{\max}, v^*)$ is a $C^2$ and increasing function satisfying

$$b(0) = 0, \quad xb(x) > 0, \quad x \neq 0 \text{ and } b'(0) > 0. \tag{4.7}$$

We consider the following bidirectional cruise controllers for $i = 1, \ldots, n$

$$F_i = \frac{1}{\beta(v_i, f_i(s))} \left( \frac{Z_i(s,v) - \mu v_{\max}^2 (v_i - f_i(s))}{v_i(v_{\max} - v_i)} + V'(s_i) - V'(s_{i+1}) \right) \tag{4.8}$$

where $\mu > 0$ is a constant,

$$\beta(v,y) := \frac{v_{\max}^3(v+y) - 2v_{\max}^2 yv}{2(v_{\max} - v)^2 v^2}, \quad v, y \in (0, v_{\max}) \tag{4.9}$$

that satisfies $\beta(v, y) > 0$ for all $v, y \in (0, v_{\max})$ and

$$Z_i(s,v) = -v_{\max}^2 b'\left(V'(s_{i+1}) - V'(s_i)\right)\left(V''(s_{i+1})(v_i - v_{i+1}) - V''(s_i)(v_{i-1} - v_i)\right) \quad i=1,...,n. \quad (4.10)$$

The construction of the Lyapunov function (4.5) is based on energy-like arguments. The first summation is a kinetic energy-like term that tends to infinity when $v_i$ tends to zero or to $v_{\max}$, while the second summation represents the potential energy. The term $-\mu(v_i - f_i(s))$ in (4.8) acts as a friction term that drives each vehicle to the spacing-dependent desired speed $f_i(s)$. Moreover, the terms $V'_{i-1}(s_i) - V'_{i+1}(s_{i+1})$ are used for collision avoidance with respect to the preceding and following vehicle, respectively.

The controller (4.8) is similar to the 2-D (lane-free) cruise controller proposed in [18], when the latter is restricted to 1-D motion. However, there are some essential differences. Here we are using measurements only from the preceding and the following vehicle while the lane-free controller in [18] requires measurements from all vehicles within the interaction distance. Another difference is the fact that here the function $b$ involved in (4.6) is an increasing function allowing desired speed $f_i(s)$ to take values greater than $v^*$. The effect of this change is essential, as the controller contains a non-zero "viscosity" term given by (4.10).

The following lemma establishes that $H(s,v)$ is a size function for (2.2), (4.8).

**Lemma 2:** *Let constants $v_{\max} > 0$, $v^* \in (0, v_{\max})$, $\lambda > L > 0$ be given. Define the function $H : \Omega \to \mathbb{R}_+$ by means of (4.5), where $\Omega$ is given by (2.3). Moreover, for every $r > \inf(H(\Omega))$ define the set*

$$S_r := \{(s,v) \in \Omega : H(s,v) \le r\}. \quad (4.11)$$

*Then for every $r \in H(\Omega)$ there exist constants $c \in (L, \lambda)$, $\underline{v} \in (0, v^*)$ and $\overline{v} \in (v^*, v_{\max})$, such that the following inequalities hold for all $(s,v) \in S_r$:*

$$R - c \ge \sum_{i=2}^{n} s_i, \quad s_i \ge c, \quad i = 2,...,n \quad (4.12)$$

$$\underline{v} \le v_i \le \overline{v}, \quad i = 1,...,n. \quad (4.13)$$

Notice now that for every $r > \inf(H(\Omega)) \ge 0$ there exists a compact set $K \subset \Omega$, namely the set

$$K = \left\{ \begin{array}{l} (s_2,...,s_n, v_1,...,v_n) \in \mathbb{R}^{2n-1}: \\ \min_{i=2,...,n}(s_i) \ge c, c + \sum_{i=2}^{n} s_i \le R \\ \max_{i=1,...,n}(v_i) \le \overline{v}, \min_{i=1,...,n}(v_i) \ge \underline{v} \end{array} \right\}$$

where $c \in (L, \lambda)$, $\underline{v} \in (0, v^*)$ and $\overline{v} \in (v^*, v_{\max})$ are provided by Lemma 2, such that $\{(s,v) \in \Omega : H(s,v) \le r\} \subseteq K$. Therefore, $H(s,v)$ is a size function for (2.2), (4.8).

Depending on the magnitude of the interaction distance $\lambda$ and the length of the ring-road $R$, we can distinguish the following two cases: (i) $R \ge \lambda n$ and (ii) $R < \lambda n$. These cases have different implications on the equilibrium points of the system (2.2), (4.8) and their stability properties.

<u>Case 1:</u> $R \ge \lambda n$ (The length of the ring-road is greater or equal to the sum of all interaction distances of the $n$ vehicles)

In this case, due to (4.2), every point in the set

$$E = \left\{ (s,v) \in \mathbb{R}^{2n-1} : \begin{array}{l} s_i \geq \lambda, i = 2,...,n, \lambda + \sum_{i=2}^{n} s_i \leq R \\ v_i = v^*, i = 1,...,n \end{array} \right\} \quad (4.14)$$

is an equilibrium point for (2.2), (4.8). Notice that $E$ is a singleton when $R = \lambda n$, namely $E = \left\{ \left( n^{-1} R 1_{n-1}, v^* 1_n \right) \right\}$.

Case 2: $R < \lambda n$ (The length of the ring-road is less than the sum of all interaction distances of the $n$ vehicles)

In this case, there is a unique equilibrium point, the point $(s,v) = \left( n^{-1} R 1_{n-1}, v^* 1_n \right)$ and

$$E = \left\{ \left( n^{-1} R 1_{n-1}, v^* 1_n \right) \right\}. \quad (4.15)$$

The set $E$ in (4.15) shows that the equilibrium inter-vehicle distance would be equal to $R/n$.

When $R \geq n\lambda$, we can use $H(s,v)$ as a Lyapunov function. Indeed, (4.1), (4.2), (4.3), (4.5), (4.6), (4.7), and (4.11) imply that $H(s,v) = 0$ when $(s,v) \in E$ and $H(s,v) > 0$ when $(s,v) \in \Omega$ with $(s,v) \notin E$. The following theorem establishes uniform global asymptotic stability of the set $E \subset \Omega$ when $R \geq \lambda n$.

**Theorem 2:** *Suppose that $R \geq n\lambda$ and consider the set $E$ given by (4.14). Then, there exist functions $\sigma \in KL$, $a \in K_\infty$ such that every solution $(s(t), v(t)) \in \Omega$ of (2.2) with (4.8) is defined for all $t \geq 0$ and satisfies the following estimate for $t \geq 0$:*

$$a\left( dist((s(t), v(t)), E) \right) \leq H(s(t), v(t)) \leq \sigma(H(s(0), v(0)), t). \quad (4.16)$$

We will consider now the case $R < \lambda n$. Define the function

$$U(s,v) = H(s,v) - nV\left( n^{-1} R \right), \quad (s,v) \in \Omega \quad (4.17)$$

where $H(s,v)$ is given by (4.5). The following proposition establishes that $U(s,v) > 0$ for $(s,v) \in \Omega \setminus E$ and $U(s,v) = 0$ for $(s,v) \in E$, when $E$ is given by (4.15).

**Proposition 1:** *Suppose that $R < \lambda n$ and consider the sets $\Omega$, $E$ given by (2.3) and (4.15), respectively. Then, the following hold for the function $U(s,v)$ defined by (4.17):*

$$U(s,v) > 0, \text{ for } (s,v) \in \Omega \setminus E \quad (4.18)$$

$$U(s,v) = 0, \text{ for } (s,v) \in E. \quad (4.19)$$

*Moreover, there exists a neighborhood $N$ of $(s,v) = \left( n^{-1} R 1_{n-1}, v^* 1_n \right)$ and constants $\alpha_1, \alpha_2, \alpha_3 > 0$ such that*

$$\alpha_1 \left| \left( s - n^{-1} R 1_{n-1}, v - v^* 1_n \right) \right|^2 \leq U(s,v) \leq \alpha_2 \left| \left( s - n^{-1} R 1_{n-1}, v - v^* 1_n \right) \right|^2, \text{ for all } (s,v) \in N \quad (4.20)$$

$$\dot{U}(s,v) \leq -\alpha_3 U(s,v), \text{ for all } (s,v) \in N \quad (4.21)$$

Proposition 1 suggests that $U(s,v)$ can be considered as a Lyapunov function for (2.2) with (4.8), while estimates (4.20) and (4.21) establish local exponential stability of $E$ defined by (4.15). Finally, notice that since $H(s,v)$ is a size function, by definition (4.17), it follows that $U(s,v)$ is also a size function.

**Theorem 3:** *Suppose that $R < n\lambda$ and consider the set $E$ is given by (4.15). Then, there exist functions $\sigma \in KL$, $a \in K_\infty$ such that every solution $(s(t),v(t)) \in \Omega$ of (2.2) with (4.8) is defined for all $t \geq 0$ and satisfies the following estimates for $t \geq 0$:*

$$a\left(\left\|\left(s(t)-n^{-1}R1_{n-1}, v(t)-v^*1_n\right)\right\|\right) \leq U(s(t),v(t)) \leq \sigma(U(s(0),v(0)),t). \tag{4.22}$$

*Moreover, there exists a size function $\kappa : \Omega \to \mathbb{R}_+$ with $\kappa(n^{-1}R1_n, v^*1_n) = 0$ and a constant $\bar{\omega} > 0$ such that every solution $(s(t),v(t)) \in \Omega$ of (2.2) with (4.8) satisfies*

$$\left|\left(s(t)-n^{-1}R1_{n-1}, v(t)-v^*1_n\right)\right| \leq \exp(-\bar{\omega}t)\kappa(s(0),v(0)), \text{ for } t \geq 0. \tag{4.23}$$

Inequalities (4.22) in Theorem 3 show that the equilibrium $(s,v) = \left(n^{-1}R1_{n-1}, v^*1_n\right)$ is uniformly globally asymptotically stable, while (4.23) indicates that every solution if (2.2) with the controller (4.8) converges exponentially to the equilibrium with rate $\bar{\omega}$. This is a direct consequence of Proposition 1, inequality (4.22), and Theorem 1. More specifically, the constant $\bar{\omega}$ in (4.23) is given by

$$\bar{\omega} = \min\left(\mu, \vartheta^2 b'(0) V''(n^{-1}R)\mu_n\right)$$

where $\mu > 0$ is the controller gain in (4.8), $\vartheta \in (0,1)$, $b(x)$ satisfying (4.7), and $\mu_n > 0$ defined by

$$\mu_n = \min\left\{\sum_{i=1}^n (x_i - x_{i-1})^2 : x = (x_1,...,x_n) \in \mathbb{R}^n, x_0 = x_n, |x| = 1, \sum_{i=1}^n x_i = 0\right\} \text{ for } n = 2,3,...$$

Notice that for every $x \in \mathbb{R}^n$ with $\sum_{i=1}^n x_i = 0$, it holds that $(x_1 - x_n)^2 + \sum_{i=2}^n (x_i - x_{i-1})^2 \geq \mu_n |x|^2$. Indeed, the left-hand side of the previous inequality forms a quadratic expression $x'Cx$ where $C$ is a circulant matrix for which $x'$ is orthogonal to $1_n$. Then, $\mu_n$ is the smallest non-zero eigenvalue of $C$, namely, $\mu_n = 2(1 - \cos(2\pi n^{-1})) > 0$, $n \geq 2$. The latter implies that the rate of convergence depends on both the parameters and functions of the controller as well as the number of vehicles.

## 5. Properties of the Cruise Controller on an Open Road

In this section we discuss some properties of the cruise controller (4.8) for the case of a string of $n$ vehicles on an open road. Notice first that since vehicles do not interact with each other when they have distance greater than a threshold $\lambda > L$ (recall (4.2)), the control problem that we study can be described as the global stabilization of the set

$$E_S = \left\{(s,v) \in \mathbb{R}^{2n-1} : \begin{matrix} s_i \geq \lambda, i = 2,...,n \\ v_i = v^*, i = 1,...,n \end{matrix}\right\} \tag{5.1}$$

where $E_S \subset \Omega_S$ is the set of equilibrium points of the closed-loop system (2.2) with (4.8), and $\Omega_S$ is given by (2.4). It should be noted that the closed set $E_S \subset \Omega_S$ is not compact. In what follows, we use the convention $s_1 = s_{n+1} = +\infty$ with $V'(s_1) = V''(s_1) = V'(s_{n+1}) = V''(s_{n+1}) = 0$.

Define the Lyapunov function

$$H_S(s,v) := \frac{v_{max}^2}{2} \sum_{i=1}^{n} \frac{(v_i - f_i(s))^2}{v_i(v_{max} - v_i)} + \sum_{i=2}^{n} V(s_i) \tag{5.2}$$

where $f(s)$ is given by (4.6) and $V(s)$ satisfies (4.1), (4.2), and (4.3).

The following result shows that $H_S(s,v)$ presents some features that characterize size functions (see [28]) for the state-space $\Omega_S$, but it is not itself a size function.

**Lemma 3:** *Let constants $v_{max} > 0$, $v^* \in (0, v_{max})$, $\lambda > L > 0$ be given. Define the function $H_S : \Omega_S \to \mathbb{R}_+$ by means of (5.2), where $\Omega_S$ is given by (2.4). Moreover, for every $r > 0$ define the set*

$$S_r := \{(s,v) \in \Omega : H_S(s,v) \leq r\}. \tag{5.3}$$

*Let $F_i$ for $i = 1,...,n$ be given by (4.8). Then for every $r > 0$ there exist constants $A > 1$, $\xi > 0$, $\underline{v} \in (0, v^*)$, and $\overline{v} \in (v^*, v_{max})$, such that the following inequalities hold for all $(s,v) \in S_r$:*

$$s_i \geq AL, \quad i = 2,...,n \tag{5.4}$$

$$\underline{v} \leq v_i \leq \overline{v}, \quad i = 1,...,n \tag{5.5}$$

$$|F_i| \leq \xi, \quad i = 1,...,n. \tag{5.6}$$

The proof of Lemma 3 is almost identical to the proof of Lemma 2 in [18] and is omitted. The following theorem guarantees that the closed-loop system (2.1) with (4.8) is well-defined and that the equilibrium set $E_S \subset \Omega_S$ is globally asymptotically stable.

**Theorem 4:** *There exist functions $\bar{\sigma} \in KL$, $\bar{a} \in K_\infty$ such that every solution $(s(t), v(t)) \in \Omega_S$ of (2.2) with (4.8) is defined for all $t \geq 0$ and satisfies the following estimate for $t \geq 0$:*

$$a\big(dist((s(t), v(t)), E_S)\big) \leq H_S(s(t), v(t)) \leq \sigma(H_S(s(0), v(0)), t). \tag{5.7}$$

**Remark: (i)** It is noted that for every $(s,v) \in \Omega_S$ the inequality

$$dist\big((s,v), E_S\big) \leq \sqrt{n(v_{max} - v^*)^2 + n(v^*)^2 + n(\lambda - L)^2}$$

holds. In other words, the whole state space is contained in a neighborhood of the equilibrium set. This implies that it is not possible to show that there exists a function $\tilde{a} \in K_\infty$ for which the inequality $H_S(s,v) \leq \tilde{a}(dist((s,v), E_S))$ holds for all $(s,v) \in \Omega_S$. Indeed, if this were true, then the inequality $dist((s,v), E_S) \leq \sqrt{n(v_{max} - v^*)^2 + n(v^*)^2 + n(\lambda - L)^2}$ would imply that $H_S(s,v)$ is bounded, which is not true.

**(ii)** Estimate (5.7) guarantees stability with respect to two measures (see [32]). However, estimate (5.7) does not guarantee exponential convergence of the solutions to the set of equilibrium points. Moreover, since the set $E_S \subset \Omega_S$ is unbounded, estimate (5.7) does not guarantee Lagrange stability (i.e., uniform boundedness of solutions for compact sets of initial data). However, we are in a position to prove the following proposition (its proof is omitted due to page limitations).

**Proposition 2:** *Every solution of (2.2) with (4.8) satisfies the following inequalities for all $t \geq 0$ and $i = 2,...,n$:*

$$s_i(t) \leq \max(s_i(0), \lambda) + \frac{v_{\max} \sqrt{2H_S(s(0), v(0))}}{2\mu \min(v^*, v_{\max} - v^*)}. \tag{5.8}$$

The following result follows directly from Lemma 1, Theorem 1 and Proposition 1.

**Corollary 1 (Uniform Lagrange Stability):** *For every compact set $S \subset \Omega_S$, there exists a compact set $K \subset \Omega_S$ such that the following implication holds for every solution of (2.2) with (4.8):*

*If $(s(0), v(0)) \in S$ then $(s(t), v(t)) \in K$ for all $t \geq 0$.*

As remarked above, Theorem 1 does not guarantee exponential convergence of the solutions to $E_S \subset \Omega_S$. However, there are certain initial conditions for which exponential convergence can be guaranteed. These are the initial conditions for which the vehicles start with sufficiently large distances between them. The following proposition shows this fact.

**Proposition 3:** *Suppose that $\min_{i=2,...,n}(s_i(0)) \geq \lambda + \frac{v_{\max}}{\mu} \Gamma$, where $\Gamma = \max_{i=1,...,n} \left( \frac{|v_i(0) - f_i(s(0))|}{\sqrt{(v_{\max} - v_i(0))v_i(0)}} \right)$. Then, the unique solution of (2.2) with (4.8) satisfies the following estimates for all $t \geq 0$:*

$$\max_{i=1,...,n}(|v_i(t) - v^*|) \leq \frac{v_{\max}}{2} \Gamma \exp(-\mu t) \tag{5.9}$$

$$\min_{i=2,...,n}(s_i(t)) \geq \lambda. \tag{5.10}$$

# 6. Proofs

**Proof of Lemma 2:** Let $r > \inf(H(\Omega)) \geq 0$ be given. Using (4.5) and the fact that $s_1 = R - \sum_{i=2}^{n} s_i$, we obtain for all $(s, v) \in S_r$:

$$V\left(R - \sum_{i=2}^{n} s_i\right) \leq r, \quad V(s_i) \leq r, \quad i = 2,...,n.$$

By virtue of (4.2), (4.3) there exists a unique $c \in (L, \lambda)$ with $V(c) = r$. Using (4.3) and the above inequalities, we obtain (4.12).

Next, from (4.2) and (4.3) we also get $V'(c) \leq V'(s_i) \leq 0$ for $i = 1,...,n$. Using (4.6), we obtain for $i = 1,...,n$:

$$0 < v^* - b(-V'(c)) \le f_i(s) \le v^* - b(V'(c)) < v_{max}. \tag{6.1}$$

Using (4.5) we also obtain for all $(s,v) \in S_r$ and $i = 1,...,n$:

$$\frac{(v_i - f_i(s))^2}{v_i(v_{max} - v_i)} \le \frac{2r}{v_{max}^2}.$$

Therefore, we get for all $(s,v) \in S_r$ and $i = 1,...,n$:

$$\frac{v_{max} f_i^2(s)}{v_{max} f_i(s) + r + \sqrt{r^2 + 2rv_{max} f_i(s) - 2rf_i^2(s)}} \le v_i$$

$$v_i \le v_{max} \frac{v_{max} f_i(s) + r + \sqrt{r^2 + 2rv_{max} f_i(s) - 2rf_i^2(s)}}{v_{max}^2 + 2r}$$

The above inequalities combined with (6.1) give the estimates (4.13) with

$$\underline{v} = \frac{v_{max}(v^* - b(-V'(\rho)))^2}{v_{max}^2 + r + \sqrt{r^2 + 2rv_{max}^2}} \text{ and }$$

$$\overline{v} = v_{max} \frac{v_{max}(v^* - b(V'(c))) + r + \sqrt{r^2 + 2r(v^* - b(V'(c)))(v_{max} - v^* + b(V'(c)))}}{v_{max}^2 + 2r}.$$

The proof is complete. ◁

**Proof of Theorem 2**: From (2.1), (4.5), and (4.6), we obtain for all $(s,v) \in \Omega$

$$\dot{H}(s,v) = \sum_{i=1}^{n} V'(s_i)(v_{i-1} - v_i) + v_{max}^2 \sum_{i=1}^{n} \frac{(v_i - f_i(s))(v_{max}(v_i + f_i(s)) - 2v_i f_i(s))}{2v_i^2(v_{max} - v_i)^2} F_i$$

$$+ v_{max}^2 \sum_{i=1}^{n} \frac{(v_i - f_i(s))}{v_i(v_{max} - v_i)} b'(V'(s_{i+1}) - V'(s_i))(V''(s_{i+1})(v_i - v_{i+1}) - V''(s_i)(v_{i-1} - v_i))$$

The above equation combined with (4.6), (4.9) and (4.10) gives

$$\dot{H}(s,v) = -\sum_{i=1}^{n} \frac{(v_i - f_i(s))}{v_i(v_{max} - v_i)} Z_i(s,v)$$
$$+ \sum_{i=1}^{n} V'(s_i)(v_{i-1} - v_i) + \sum_{i=1}^{n} (v_i - f_i(s)) \beta(v_i, f_i(s)) F_i \tag{6.2}$$

Using now (6.2), (4.6), and (4.8) we have that

$$\dot{H}(s,v) = -\mu v_{max}^2 \sum_{i=1}^{n} \frac{(v_i - f_i(s))^2}{v_i(v_{max} - v_i)}$$
$$+ \sum_{i=1}^{n} V'(s_i)(v_{i-1} - v^*) - \sum_{i=1}^{n} V'(s_{i+1})(v_i - v^*) \tag{6.3}$$
$$- \sum_{i=1}^{n} b(V'(s_{i+1}) - V'(s_i))(V'(s_{i+1}) - V'(s_i))$$

Since $v_0 = v_n$ and $s_{n+1} = s_1$, we get $\sum_{i=1}^{n} V'(s_i)(v_{i-1} - v^*) = \sum_{i=}^{n} V'(s_{i+1})(v_i - v^*)$, which, in conjunction with (6.3) gives

$$\dot{H}(s,v) = -\mu v_{max}^2 \sum_{i=1}^{n} \frac{(v_i - f_i(s))^2}{v_i(v_{max} - v_i)} - \sum_{i=1}^{n} b(V'(s_{i+1}) - V'(s_i))(V'(s_{i+1}) - V'(s_i)). \quad (6.4)$$

It follows from (4.7) and (6.4) that $\dot{H}(s,v) < 0$ for all $(s,v) \in \Omega \setminus E$. Notice that (6.4) implies that $H(s(t), v(t)) \le H(s(0), v(0))$ for as long as the solution is defined. By virtue of Lemma 2 and Theorem 1 in [18], we conclude that every solution of (2.2), (4.8) is defined for all $t \ge 0$ and satisfies $(s(t), v(t)) \in \Omega$.

Applying Lemma 1 with $O = \Omega$, $D = E$, $W(s,v) = -\dot{H}(s,v)$ and $Q(s,v) = H(s,v)$ it follows that there exists a positive definite and globally Lipschitz function $\rho: \mathbb{R}_+ \to \mathbb{R}_+$ such that

$$\dot{H}(s,v) \le -\rho(H(s,v)) \text{ for all } (s,v) \in \Omega. \quad (6.5)$$

Thus, for every solution $(s(t), v(t)) \in \Omega$, $t \ge 0$, inequality (6.5) and Lemma 4.4 in [20] imply the existence of $\sigma \in KL$ such that

$$H(s(t), v(t)) \le \sigma(H(s(0), v(0)), t) \text{ for all } t \ge 0. \quad (6.6)$$

Since $H(s,v)$ is a size function and $H(\Omega) = \mathbb{R}_+$, applying again Lemma 1 with $O = \Omega$, $D = E$, $W(s,v) = dist((s,v), E)$, and $Q(s,v) = H(s,v)$, we conclude that there exists $\zeta \in K_\infty$ such that the following inequalities hold for all $(s,v) \in \Omega$:

$$dist((s,v), E) \le \zeta(H(s,v)). \quad (6.7)$$

Inequality (4.16) follows from (6.6) and (6.7) with $a = \zeta^{-1}$. This concludes the proof. $\square$

**Proof of Proposition 1:** Notice first that (4.5), (4.6), (4.7), and definition (4.17) imply that $U(s,v) = 0$ when $(s,v) \in E = \{(n^{-1}R1_{n-1}, v^*1_n)\}$. Next, define the function

$$G(s) = \sum_{i=1}^{n} V(s_i) - nV(n^{-1}R), \text{ for all } s \in \left\{ (y_1, ..., y_n) \in \mathbb{R}^n : \min_{i=1,...,n}(y_i) > L, \sum_{i=1}^{n} y_i = R \right\} \quad (6.8)$$

where $V(s)$ satisfies (4.1), (4.2), and (4.3). Since

$$V(s_i) = V(n^{-1}R) + V'(n^{-1}R)(s_i - n^{-1}R) + (s_i - n^{-1}R)^2 \int_0^1 \int_0^z V''(n^{-1}R + r(s_i - n^{-1}R)) dr dz \quad (6.9)$$

for all $s \in \left\{ (y_1, ..., y_n) \in \mathbb{R}^n : \min_{i=1,...,n}(y_i) > L, \sum_{i=1}^{n} y_i = R \right\}$, we get:

$$G(s) = V'(n^{-1}R)\left(\sum_{i=1}^{n} s_i - R\right) + \sum_{i=1}^{n} (s_i - n^{-1}R)^2 \int_0^1 \int_0^z V''(n^{-1}R + r(s_i - n^{-1}R)) dr dz$$
$$= \sum_{i=1}^{n} (s_i - n^{-1}R)^2 \int_0^1 \int_0^z V''(n^{-1}R + r(s_i - n^{-1}R)) dr dz \quad (6.10)$$

The above equation in conjunction with (4.3) implies that $Q(s) > 0$ for all $s \in \left\{ (y_1,...,y_n) \in \mathbb{R}^n : \min_{i=1,...,n}(y_i) > L, \sum_{i=1}^{n} y_i = R \right\}$ with $s \neq n^{-1}R1_n$. Therefore, it follows that

$$V\left(R - \sum_{i=2}^{n} s_i\right) + \sum_{i=2}^{n} V(s_i) > nV\left(n^{-1}R\right), \text{ for all}$$

$$s \neq n^{-1}R1_{n-1}, s \in \left\{ (y_2,...,y_n) \in \mathbb{R}^{n-1} : \min_{i=2,...,n}(y_i) > L, L + \sum_{i=1}^{n} y_i < R \right\}.$$

Consequently, definition (4.5) implies that $U(s,v) = H(s,v) - nV\left(n^{-1}R\right) > 0$ when $(s,v) \in \Omega$ with $(s,v) \neq \left(n^{-1}R1_{n-1}, v^*1_n\right)$.

Next, we show that estimate (4.20) holds in a neighborhood of $(s,v) = \left(n^{-1}R1_{n-1}, v^*1_n\right)$. From definition (6.10), we have that for every $\vartheta \in (0,1)$, there exists $\delta > 0$ such that for every $s_i$ with $\min_{i=1,...,n}(s_i) > L$, $\sum_{i=1}^{n} s_i = R$, and

$$\sum_{i=1}^{n} \left(s_i - n^{-1}R\right)^2 \leq \delta \tag{6.11}$$

it holds that

$$\vartheta \frac{V''\left(n^{-1}R\right)}{2} \sum_{i=1}^{n} \left(s_i - n^{-1}R\right)^2 \leq G(s) \leq \frac{V''\left(n^{-1}R\right)}{2\vartheta} \sum_{i=1}^{n} \left(s_i - n^{-1}R\right)^2. \tag{6.12}$$

Notice that by definitions (4.5), (4.17) and (6.10), we have that

$$U(s,v) = \frac{v_{\max}^2}{2} \sum_{i=1}^{n} \frac{(v_i - f_i(s))^2}{v_i(v_{\max} - v_i)} + G(s). \tag{6.13}$$

Using (6.12), inequality $v_i(v_{\max} - v_i) \leq \frac{v_{\max}^2}{4}$ that holds for all $v_i \in (0, v_{\max})$, (4.6), (6.12), and (6.13), we have that

$$\begin{aligned} U(s,v) &\geq 2\sum_{i=1}^{n}\left(v_i - v^* + b(V'(s_{i+1}) - V'(s_i))\right)^2 + \vartheta \frac{V''\left(n^{-1}R\right)}{2} \sum_{i=1}^{n}\left(s_i - n^{-1}R\right)^2 \\ &= 2\sum_{i=1}^{n}(v_i - v^*)^2 + 4\sum_{i=1}^{n}(v_i - v^*)b(V'(s_{i+1}) - V'(s_i)) \\ &\quad + 2\sum_{i=1}^{n}\left(b(V'(s_{i+1}) - V'(s_i))\right)^2 + \vartheta \frac{V''\left(n^{-1}R\right)}{2} \sum_{i=1}^{n}\left(s_i - n^{-1}R\right)^2 \end{aligned} \tag{6.14}$$

Let $\varepsilon > 1$. By adding and subtracting terms, using (6.14) and the inequality $\frac{1}{\varepsilon}\sum_{i=1}^{n}\left(v_i - v^* + b(V'(s_{i+1}) - V'(s_i))\right)^2 \geq 0$, it follows that the following estimate holds

$$U(s,v) \geq 2(1-\varepsilon^{-1})\sum_{i=1}^{n}(v_i - v^*)^2 + \vartheta \frac{V''(n^{-1}R)}{2}\sum_{i=1}^{n}(s_i - n^{-1}R)^2 \qquad (6.15)$$
$$-2(\varepsilon-1)\sum_{i=1}^{n}\left(b(V'(s_{i+1})-V'(s_i))\right)^2$$

Notice now that we can rewrite the term $b(V'(s_{i+1})-V'(s_i))$ as follows

$$b(V'(s_{i+1})-V'(s_i)) = \int_{s_i}^{s_{i+1}} b'(V'(y)-V'(s_i))V''(y)dy.$$

The latter, in conjunction with (4.3, (4.7), and (6.9), gives that

$$\left|b(V'(s_{i+1})-V'(s_i))\right| \leq \vartheta^{-1}b'(0)V''(n^{-1}R)|s_{i+1}-s_i| \qquad (6.16)$$

for $\vartheta \in (0,1)$ and $\delta > 0$ as above. It follows by (6.15), (6.16) and inequality $\sum_{i=1}^{n}(s_{i+1}-s_i)^2 \leq 4\sum_{i=1}^{n}(s_i - n^{-1}R)^2$ that

$$U(s,v) \geq 2(1-\varepsilon^{-1})|v-v^*1_n|^2$$
$$+ \frac{V''(n^{-1}R)}{2}\vartheta\left(1-16(\varepsilon-1)\vartheta^{-3}(b'(0))^2 V''(n^{-1}R)\right)\sum_{i=1}^{n}(s_i - n^{-1}R)^2 \qquad (6.17)$$

establishing the first inequality of (4.20) with

$$\alpha_1 = \min\left(2(1-\varepsilon^{-1}), \frac{V''(n^{-1}R)}{2}\vartheta\left(1-16(\varepsilon-1)\vartheta^{-3}(b'(0))^2 V''(n^{-1}R)\right)\right).$$

The existence of a constant $\alpha_2$ that satisfies the right-hand side inequality in (4.20) in a neighborhood of $E$ for which (6.11) holds, is a direct consequence of (4.18), (4.19), and the fact that $U \in C^3((L,\lambda); \mathbb{R}_+)$.

We proceed now with the proof of (4.23). Define

$$g(x) = xb(x), \quad x \in \mathbb{R}. \qquad (6.18)$$

Due to (4.7) and the fact that $g'(x) = b(x) + xb'(x)$, and $g''(x) = 2b'(x) + xb''(x)$, we have for $\vartheta \in (0,1)$ as above and $\delta > 0$ satisfying (6.11) that

$$\vartheta b'(0)(V'(s_{i+1})-V'(s_i))^2 \leq g(V'(s_{i+1})-V'(s_i)) \leq \vartheta^{-1}b'(0)(V'(s_{i+1})-V'(s_i))^2$$

and consequently, by (6.9) we have

$$\vartheta b'(0)(V''(n^{-1}R))^2 (s_{i+1}-s_i)^2 \leq g(V'(s_{i+1})-V'(s_i)) \leq \vartheta^{-1}b'(0)(V''(n^{-1}R))^2 (s_{i+1}-s_i)^2. \qquad (6.19)$$

Define for $n = 2,3,...$

$$\mu_n = \min\left\{\sum_{i=1}^{n}(x_i - x_{i-1})^2 : x = (x_1,...,x_n) \in \mathbb{R}^n, x_0 = x_n, |x|=1, \sum_{i=1}^{n}x_i = 0\right\} \qquad (6.20)$$

and notice that $\mu_n > 0$ for all $n = 2,3,...$ and that for every $x \in \mathbb{R}^n$ with $\sum_{i=1}^{n}x_i = 0$ it holds that

$$(x_1 - x_n)^2 + \sum_{i=2}^{n}(x_i - x_{i-1})^2 \geq \mu_n |x|^2. \tag{6.21}$$

Since, $\sum_{i=1}^{n}(s_i - n^{-1}R) = 0$, using (6.21) with $x_i = s_i - n^{-1}R$ it holds that

$$(s_1 - s_n)^2 + \sum_{i=2}^{n}(s_i - s_{i-1})^2 \geq \mu_n \sum_{i=1}^{n}(s_i - n^{-1}R)^2. \tag{6.22}$$

Define

$$L(s) = \sum_{i=1}^{n} g\big(V'(s_{i+1}) - V'(s_i)\big), \text{ for } s \in \left\{(y_2,...,y_n) \in \mathbb{R}^{n-1}: \min_{i=2,...,n}(y_i) > L, L + \sum_{i=1}^{n} y_i < R\right\}.$$

It follows from (6.19) and (6.22) that

$$\vartheta b'(0)\big(V''(n^{-1}R)\big)^2 \mu_n \sum_{i=1}^{n}(s_i - n^{-1}R)^2 \leq L(s) \leq 4\vartheta^{-1} b'(0)\big(V''(n^{-1}R)\big)^2 \sum_{i=1}^{n}(s_i - n^{-1}R)^2. \tag{6.23}$$

Using definitions (4.5) and (4.17), and (6.4), it follows that

$$\dot{U}(s,v) = \dot{H}(s,v) = -\mu v_{\max}^2 \sum_{i=1}^{n} \frac{(v_i - f_i(s))^2}{v_i(v_{\max} - v_i)} - \sum_{i=1}^{n} b\big(V'(s_{i+1}) - V'(s_i)\big)\big(V'(s_{i+1}) - V'(s_i)\big) \tag{6.24}$$

and consequently, from (6.12), (6.13), and (6.23) we get

$$\begin{aligned}
\dot{U}(s,v) &\leq -\mu v_{\max}^2 \sum_{i=1}^{n} \frac{(v_i - f_i(s))^2}{v_i(v_{\max} - v_i)} - \vartheta b'(0)\big(V''(n^{-1}R)\big)^2 \mu_n \sum_{i=1}^{n}(s_i - n^{-1}R)^2 \\
&= -2\mu \frac{v_{\max}^2}{2} \sum_{i=1}^{n} \frac{(v_i - f_i(s))^2}{v_i(v_{\max} - v_i)} - 2\vartheta^2 b'(0) V''(n^{-1}R) \mu_n \frac{V''(n^{-1}R)}{2\vartheta} \sum_{i=1}^{n}(s_i - n^{-1}R)^2 \\
&\leq -2\min\big(\mu, \vartheta^2 b'(0) V''(n^{-1}R) \mu_n\big) \left(\frac{v_{\max}^2}{2} \sum_{i=1}^{n} \frac{(v_i - f_i(s))^2}{v_i(v_{\max} - v_i)} + \frac{V''(n^{-1}R)}{2\vartheta} \sum_{i=1}^{n}(s_i - n^{-1}R)^2\right) \\
&\leq -2\min\big(\mu, \vartheta^2 b'(0) V''(n^{-1}R) \mu_n\big) U(s,v)
\end{aligned} \tag{6.25}$$

Thus, (4.21) holds with $\alpha_3 = 2\min\big(\mu, \vartheta^2 b'(0) V''(n^{-1}R) \mu_n\big)$, in a neighborhood of $E$ for which (6.11) holds. This completes the proof. □

**Proof of Theorem 3:** From (6.24) and by applying Lemma 1 with $O = \Omega$, $D = E$, $W(s,v) = -\dot{U}(s,v) = -\dot{H}(s,v)$, and $Q(s,v) = U(s,v)$, it follows that there exists a positive definite and globally Lipschitz function $\rho: \mathbb{R}_+ \to \mathbb{R}_+$ such that

$$\dot{U}(s,v) \leq -\rho(U(s,v)), \text{ for all } (s,v) \in \Omega. \tag{6.26}$$

Since $U(s,v)$ is a size function and $U(\Omega) = \mathbb{R}_+$, applying again Lemma 1 with $O = \Omega$, $D = E$, $W(s,v) = dist\big((s,v), E\big) = \big|(s - n^{-1}R1_{n-1}, v - v^*1_n)\big|$, and $Q(s,v) = U(s,v)$, we conclude that there exists $\zeta \in K_\infty$ such that the following inequalities hold for all $(s,v) \in \Omega$:

$$\big|(s - n^{-1}R1_{n-1}, v - v^*1_n)\big| \leq \zeta(U(s,v)). \tag{6.27}$$

Thus, for every solution $(s(t),v(t)) \in \Omega$, $t \geq 0$, inequality (6.26) and Lemma 4.4 in [20], imply the existence of $\sigma \in KL$ such that

$$U(s(t),v(t)) \leq \sigma(U(s(0),v(0)),t) \text{ for all } t \geq 0$$

which, together with (6.27), establishes inequality (4.22) with $a = \zeta^{-1}$. Inequality (4.23) is a direct consequence of (4.22) established earlier, Proposition 1 and Theorem 1 with $\bar{\omega} = \min\left(\mu, \vartheta^2 b'(0) V''(n^{-1}R)\mu_n\right)$ for $\vartheta \in (0,1)$ and $\mu_n$ given by (6.20). This concludes the proof. □

**Proof of Theorem 4:** From (2.2), (4.5), (4.6), (4.9), (4.10), and (5.2), we obtain for all $(s,v) \in \Omega_S$

$$\dot{H}_S(s,v) = \sum_{i=1}^{n} \beta(v_i, f_i(s))(v_i - f_i(s))F_i$$

$$-\sum_{i=1}^{n}(v_i - f_i(s))\frac{Z_i(s,v)}{v_i(v_{\max} - v_i)} + \sum_{i=2}^{n} V'(s_i)(v_{i-1} - v_i)$$

Using (4.8) it follows that

$$\dot{H}_S(s,v) = -\sum_{i=1}^{n} \frac{\mu v_{\max}^2 (v_i - f_i(s))^2}{v_i(v_{\max} - v_i)} - \sum_{i=1}^{n} b(V'(s_{i+1}) - V'(s_i))(V'(s_{i+1}) - V'(s_i)). \tag{6.28}$$

It follows from (4.7) and (6.28) that $\dot{H}_S(s,v) < 0$ for all $(s,v) \in \Omega_S \setminus E_S$ with $E_S$ defined by (5.1). Notice that (6.28) implies that $H_S(s(t),v(t)) \leq H_S(s(0),v(0))$ for as long as the solution is defined. By virtue of Lemma 1 and Theorem 1 in [18] we conclude that every solution of (2.2), (4.8) is defined for all $t \geq 0$ and satisfies $(s(t),v(t)) \in \Omega$.

Let $P:(L,+\infty)^{n-1} \to (L,\lambda]^{n-1}$ be the projection mapping on $(L,\lambda]^{n-1}$. Notice that (4.5), (4.6), (6.28) and (4.2) imply that $\dot{H}_S(s,v) = \dot{H}_S(P(s),v)$, $H_S(s,v) = H_S(P(s),v)$, for all $(s,v) \in \Omega_S$. Define

$$\Lambda := \left\{ \begin{array}{l} (s_2,...,s_n,v_1,...,v_n) \in \mathbb{R}^{2n-1}: \\ L < \min_{i=2,...,n}(s_i), \max_{i=2,...,n}(s_i) \leq \lambda, \\ 0 < \min_{i=1,...,n}(v_i), \max_{i=1,...,n}(v_i) < v_{\max} \end{array} \right\}, \quad \bar{D} := \{(\lambda \mathbf{1}_{n-1}, v^* \mathbf{1}_n)\}.$$

Applying Lemma 1 with $Q(s,v) = H_S(s,v)$, $W(s,v) = -\dot{H}_S(s,v)$ and $O = \Lambda$, $D = \bar{D}$ as above, we conclude that there exists a continuous and positive definite function $\rho: \mathbb{R}_+ \to \mathbb{R}_+$ such that

$$\dot{H}_S(s,v) \leq -\rho(H_S(s,v)), \text{ for all } (s,v) \in \Lambda.$$

Since $(P(s),v) \in \Lambda$ for all $(s,v) \in \Omega_S$ we get that:

$$\dot{H}_S(s,v) \leq -\rho(H_S(s,v)), \text{ for all } (s,v) \in \Omega_S. \tag{6.29}$$

Thus, for every solution $(s(t),v(t)) \in \Omega_S$, $t \geq 0$, inequality (6.29) and Lemma 4.4 in [20], imply the existence of $\sigma \in KL$ such that

$$H_S(s(t),v(t)) \leq \sigma(H_S(s(0),v(0)),t) \text{ for all } t \geq 0. \tag{6.30}$$

We finally establish the existence of a function $\bar{a} \in K_\infty$ so that $\bar{a}(dist((s,v),E_S)) \leq H_S(s,v)$ for all $(s,v) \in \Omega_S$. Let $\tilde{P}:\mathbb{R}^{n-1} \to [\lambda,+\infty)^{n-1}$ be the projection mapping on $[\lambda,+\infty)^{n-1}$ defined by

$$\mathbb{R}^{n-1} \ni (s_2,...,s_n) \to y = (y_2,...,y_n) = \tilde{P}(s)$$

$$\text{with } y_i = \max(\lambda, s_i), \ i = 2,\ldots,n \tag{6.31}$$

and notice that for every $(s,v) \in \Omega_S$ it holds that $(\tilde{P}(s), v^* \mathbf{1}_n) \in E$. Using (6.31) and (5.1), we obtain for all $(s,v) \in \Omega_S$:

$$dist((s,v), E_S) \leq \left( |\tilde{P}(s) - s|^2 + |v - v^* \mathbf{1}_n|^2 \right)^{1/2}$$
$$\leq \left( \sum_{i=1}^{n} |v_i - v^*|^2 + \sum_{i=2}^{n} (\max(\lambda - s_i, 0))^2 \right)^{1/2} \tag{6.32}$$

Notice also that due to (4.5) and (4.6) we have

$$|v_i - v^*|^2 \leq 2|v_i - f_i(s)|^2 + 2|f_i(s) - v^*|^2$$
$$\leq 4v_{\max}^{-2} H_S(s,v) + 2B(|V'(s_{i+1}) - V'(s_i)|) \tag{6.33}$$

where $B(r) := \max\{b^2(l) : |l| \leq r\}$ for $r \geq 0$ is a non-decreasing, non-negative function which is continuous at $r = 0$ with $B(0) = 0$ (recall (4.7)). Recall that by virtue of Lemma 6.1 in [16], there exists a decreasing, continuous function $\tilde{\rho}: \mathbb{R}_+ \to (L, \lambda]$ such that

$$s_i \geq \tilde{\rho}(H(s,v)) \text{ for all } (s,v) \in \Omega_S, \ i = 2,\ldots,n \tag{6.34}$$

with $\tilde{\rho}(0) = \lambda$ and $\lim_{l \to +\infty} (\tilde{\rho}(l)) = L$. Define

$$\gamma(r) := \sup\{|V'(l)| : l \geq r\} \text{ for } r > L. \tag{6.35}$$

Using (4.2) we can guarantee that $\gamma$ defined by (6.35) is a well-defined, non-negative, non-increasing function which is continuous at $r = \lambda$ with $\gamma(r) = 0$ for $r \geq \lambda$. Combining (6.34) and (6.35), we have

$$|V'(s_i) - V'(s_{i+1})| \leq \gamma(s_i) + \gamma(s_{i+1}) \leq 2\gamma(\tilde{\rho}(H_S(s,v))). \tag{6.36}$$

Note that $\tilde{Q}(r) := 4v_{\max}^{-2} r + 2B(2\gamma(\tilde{\rho}(r)))$ for $r \geq 0$ is an increasing, non-negative function which is continuous at $r = 0$ with $\tilde{Q}(0) = 0$. By virtue of Lemma 2.4 on page 65 in [14] there exists $\bar{Q} \in K_\infty$ such that $\bar{Q}(r) \geq \tilde{Q}(r)$ for all $r \geq 0$. Combining (6.33) and (6.36) we obtain

$$|v_i - v^*|^2 \leq \bar{Q}(H_S(s,v)). \tag{6.37}$$

The function $\bar{b}(l) := \lambda - \tilde{\rho}(l)$ defined on $\mathbb{R}_+$ is an increasing, continuous function with $\bar{b}(0) = 0$ and $\lim_{l \to +\infty}(\bar{b}(l)) = \lambda - L$. Moreover, we have (due to (6.34))

$$\max(\lambda - s_i, 0) \leq \lambda - \tilde{\rho}(H_S(s,v)) = \bar{b}(H_S(s,v))$$

for all $(s,v) \in \Omega_S$ and $i = 2,\ldots,n$. Combining the above inequality, (6.37), (6.33), and (6.32) we obtain for all $(s,v) \in \Omega_S$:

$$dist((s,v), S) \leq \left( n\bar{Q}(H_S(s,v)) + (n-1)(\bar{b}(H_S(s,v)))^2 \right)^{1/2} \tag{6.38}$$

The previous definitions imply that the function $h(l) = \left(n\bar{Q}(l) + (n-1)(\bar{b}(l))^2\right)^{1/2}$ for $l \geq 0$, is of class $K_\infty$ and therefore it holds that

$$\bar{a}(dist((s,v),S)) \leq H_S(s,v) \text{ for all } (s,v) \in \Omega_S$$

with $\bar{a} = h^{-1}$. The proof is complete. ◁

**Proof of Proposition 2:** Define the function

$$p(x) = \frac{x - v^*}{\sqrt{(v_{max} - x)x}}, \text{ for } x \in (0, v_{max}). \tag{6.39}$$

The function $p : (0, v_{max}) \to \mathbb{R}$ is a $C^1$ increasing function with

$$p'(x) = \frac{v_{max} x + v^* v_{max} - 2v^* x}{2(v_{max} - x)x\sqrt{(v_{max} - x)x}} \geq \frac{4\min(v^*, v_{max} - v^*)}{v_{max}^2} > 0 \text{ for all } x \in (0, v_{max}) \text{ and}$$

$p((0, v_{max})) = \mathbb{R}$. Therefore, the inverse function of $p$ is well-defined and is a $C^1$ increasing function $p^{-1}: \mathbb{R} \to (0, v_{max})$ with $(p^{-1})'(x) \leq \frac{v_{max}^2}{4\min(v^*, v_{max} - v^*)}$ for all $x \in \mathbb{R}$. Consequently, we get for all $x, y \in \mathbb{R}$:

$$p^{-1}(x) - p^{-1}(y) \leq \frac{v_{max}^2}{4\min(v^*, v_{max} - v^*)} \max(0, x - y). \tag{6.40}$$

Define for $i = 1, \ldots, n$:

$$z_i = \frac{v_i - f_i(s)}{\sqrt{(v_{max} - v_i)v_i}}. \tag{6.41}$$

Using (4.5) and definition (6.41), we get for all $i = 1, \ldots, n$ and $(s,v) \in \Omega_S$:

$$|z_i| \leq \frac{1}{v_{max}} \sqrt{2H_S(s,v)}. \tag{6.42}$$

Suppose that $s_i \geq \lambda$ for some $i = 2, \ldots, n$. Using (4.2), (4.3), (4.6), (4.7) and the fact that $b$ is increasing, we get:

$$f_i(s) = v^* - b(V'(s_{i+1})) \geq v^*, \text{ if } i < n$$

$$f_i(s) = v^*, \text{ if } i = n$$

$$f_{i-1}(s) = v^* - b(-V'(s_{i-1})) \leq v^*, \text{ if } i > 2$$

$$f_{i-1}(s) = v^*, \text{ if } i = 2.$$

Consequently, we get from (6.41) and (6.39) when $s_i \geq \lambda$ for some $i = 2, \ldots, n$:

$$z_i \leq p(v_i) \text{ and } z_{i-1} \geq p(v_{i-1}). \tag{6.43}$$

Using (6.41), (2.2), (4.6), (4.8), (4.10) and (4.11), we get:

$$\dot{z}_1 = -\mu z_1 - \sqrt{(v_{max} - v_1)v_1} \frac{V'(s_2)}{v_{max}^2}$$

$$\dot{z}_i = -\mu z_i + \sqrt{(v_{max} - v_i)v_i}\, \frac{V'(s_i) - V'(s_{i+1})}{v_{max}^2}$$

$$\dot{z}_n = -\mu z_n + \sqrt{(v_{max} - v_n)v_n}\, \frac{V'(s_n)}{v_{max}^2}.$$

Consequently, we get from (4.3) and (4.2) when $s_i \geq \lambda$ for some $i = 2,\ldots,n$:

$$\dot{z}_i \geq -\mu z_i \text{ and } \dot{z}_{i-1} \leq -\mu z_{i-1}. \tag{6.44}$$

Using (2.2), (6.40) and (6.43), we get when $s_i \geq \lambda$ for some $i = 2,\ldots,n$:

$$\dot{s}_i = v_{i-1} - v_i \leq p^{-1}(z_{i-1}) - p^{-1}(z_i)$$

$$\leq \frac{v_{max}^2}{4\min(v^*, v_{max} - v^*)} \max(0, z_{i-1} - z_i) \tag{6.45}$$

First notice that (5.8) holds for $t = 0$. The proof of (5.8) for $t > 0$ is made by contradiction. Suppose that there exists a solution of (2.2) with (4.8) for which the following inequality holds for certain $t > 0$ and $i = 2,\ldots,n$:

$$s_i(t) > \max(s_i(0), \lambda) + \frac{v_{max}\sqrt{2H_S(s(0), v(0))}}{2\mu \min(v^*, v_{max} - v^*)}. \tag{6.46}$$

If the set $\{\tau \in [0,t]: s_i(\tau) \leq \lambda\}$ is empty then we set $T = 0$. If the set $\{\tau \in [0,t]: s_i(\tau) \leq \lambda\}$ is non-empty then we set $T = \sup(\{\tau \in [0,t]: s_i(\tau) \leq \lambda\})$. In any case, we have $T < t$ and $s_i(\tau) \geq \lambda$ for $\tau \in [T,t]$. It follows from (6.44) that

$$z_i(\tau) \geq e^{-\mu(\tau - T)} z_i(T),\ z_{i-1}(\tau) \leq e^{-\mu(\tau - T)} z_{i-1}(T),\ \text{for } \tau \in [T,t]. \tag{6.47}$$

Consequently, we get from (6.47):

$$\max(0, z_{i-1}(\tau) - z_i(\tau)) \leq e^{-\mu(\tau - T)} \max(0, z_{i-1}(T) - z_i(T)),\ \text{for } \tau \in [T,t]. \tag{6.48}$$

Using (6.45) and (6.48) we get for $\tau \in [T,t]$:

$$\dot{s}_i(\tau) \leq \frac{v_{max}^2}{4\min(v^*, v_{max} - v^*)} e^{-\mu(\tau - T)} \max(0, z_{i-1}(T) - z_i(T)).$$

Integrating the above differential inequality, we get for $\tau \in [T,t]$:

$$s_i(\tau) \leq s_i(T) + \frac{v_{max}^2}{4\mu \min(v^*, v_{max} - v^*)} \max(0, z_{i-1}(T) - z_i(T)). \tag{6.49}$$

Since $s_i(T) \leq \max(s_i(0), \lambda)$ (distinguish the cases that the set $\{\tau \in [0,t]: s_i(\tau) \leq \lambda\}$ is empty or not), we get from (6.49) for $\tau \in [T,t]$:

$$s_i(\tau) \leq \max(s_i(0), \lambda) + \frac{v_{max}^2}{4\mu \min(v^*, v_{max} - v^*)} (|z_{i-1}(T)| + |z_i(T)|). \tag{6.50}$$

Using (6.42) and (6.50) we get from (6.49) for $\tau \in [T,t]$:

$$s_i(\tau) \le \max\left(s_i(0), \lambda\right) + \frac{v_{\max}\sqrt{2H_S(s(T), v(T))}}{2\mu \min\left(v^*, v_{\max} - v^*\right)}. \tag{6.51}$$

Inequality (6.33) implies that $H_S(s(T), v(T)) \le H_S(s(0), v(0))$. Hence, we obtain from (6.51) for $\tau \in [T, t]$:

$$s_i(\tau) \le \max\left(s_i(0), \lambda\right) + \frac{v_{\max}\sqrt{2H_S(s(0), v(0))}}{2\mu \min\left(v^*, v_{\max} - v^*\right)}. \tag{6.52}$$

Estimate (6.52) contradicts inequality (6.46). The proof is complete. ◁

## 7. Simulations

In this section we demonstrate the properties of the cruise controller (4.8) for both the cases of the ring-road and the open road.

### *5.1 Application on a Ring-Road*

We consider $n = 4$ vehicles on a ring-road of length $R = 130$, with the cruise controller given by (4.8) with $L = 5$, $\mu = 0.1$, $v_{\max} = 35$, $v^* = 30$ and initial conditions $s(0) = (38, 33, 32, 27)$ and $v(0) = (31, 28, 27, 30)$. The potential $V$ satisfying (4.1), (4.2), and (4.3) is given by

$$V(s) = \begin{cases} q\dfrac{(\lambda - s)^4}{(s - L)^2} & L < s \le \lambda \\ 0 & s > \lambda \end{cases} \tag{7.1}$$

with $q = 0.1$ and $\lambda > L > 0$. We also select $b(x)$ defined by

$$b(x) = v^* + \frac{v_{\max}}{2}\left(\tanh(x + c) - 1\right) \text{ with } c = \operatorname{arctanh}\left(1 - \frac{2v^*}{v_{\max}}\right). \tag{7.2}$$

First, we demonstrate the results of Theorem 2, with $\lambda = 30$ which satisfies $R > n\lambda$. Figure 1 shows the convergence of the speeds $v_i$, $i = 1, \ldots, 4$, to the desired speed $v^*$. Figure 2 shows the evolution of the intervehicle distances and their convergence to values great or equal to $\lambda$ (recall the set of equilibrium points $E$ given by (4.14)). Finally, Figure 3 shows the convergence of the accelerations $F_i(t)$ for each vehicle.

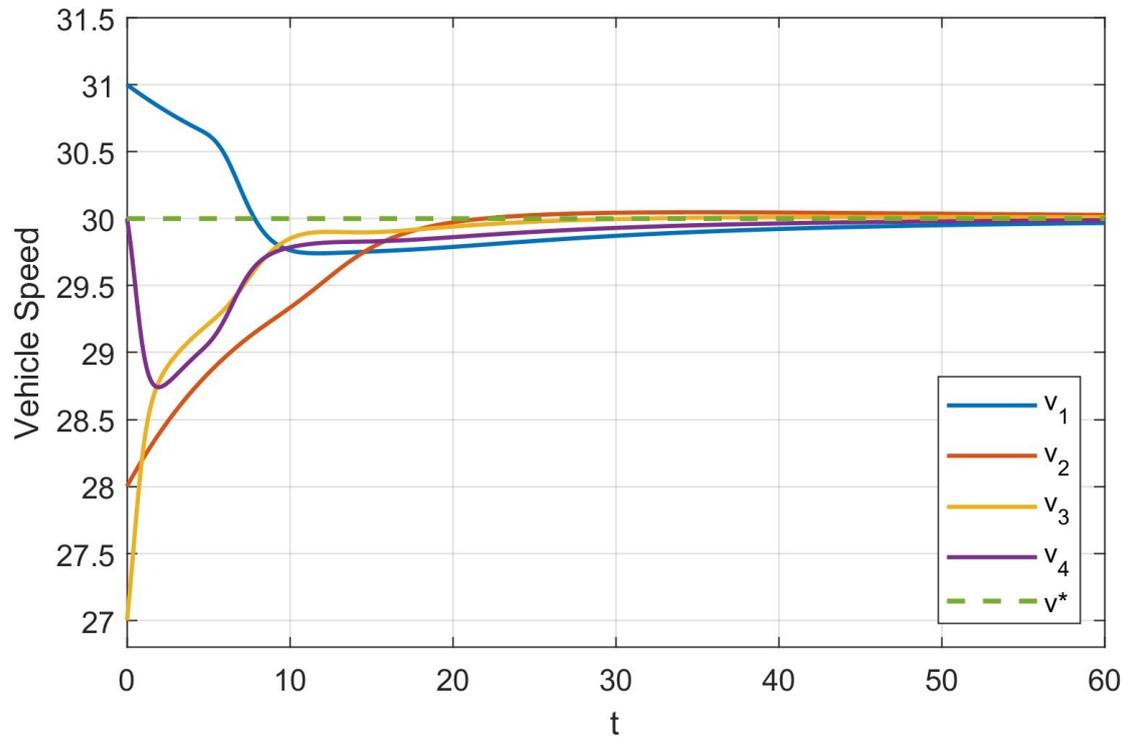

**Figure 1:** Convergence of vehicle speeds to the desired speed $v^*$ for the case $R > n\lambda$.

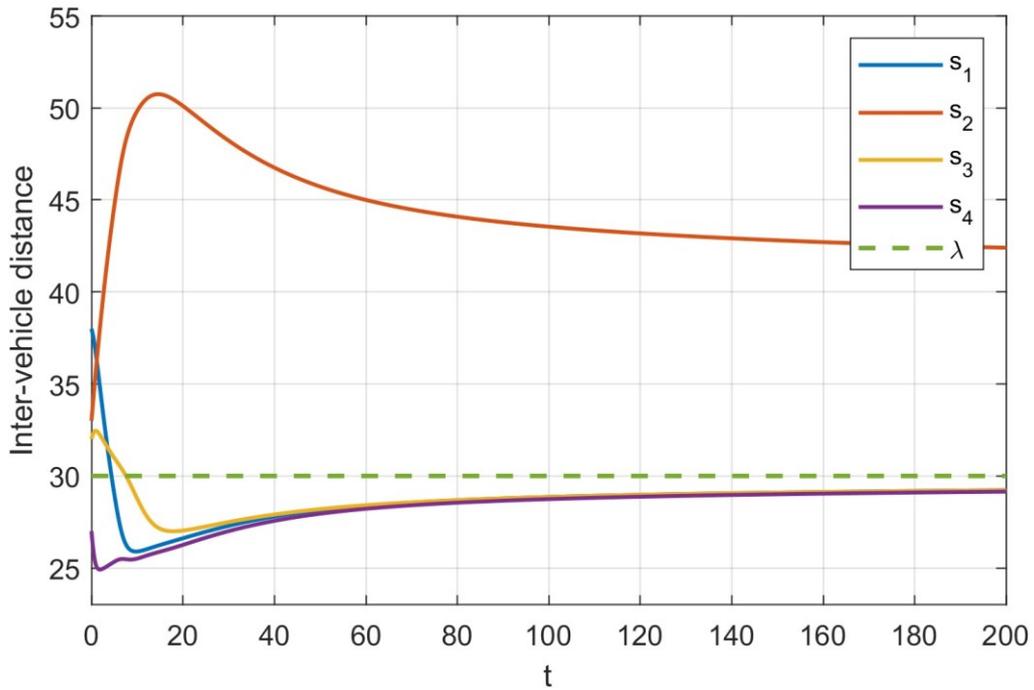

**Figure 2:** Evolution of inter-vehicle distance $s_i$ for the case $R > n\lambda$.

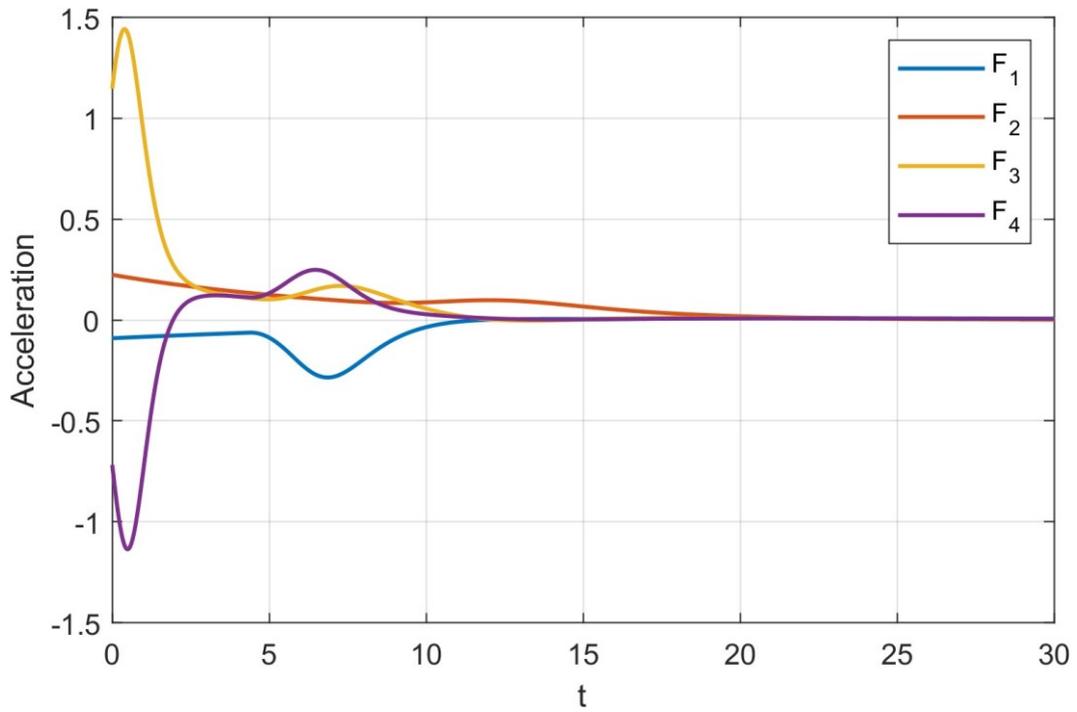

**Figure 3:** Convergence of vehicle accelerations $F_i$ for the case $R > n\lambda$.

We consider now the case where $\lambda = 40$ which satisfies $R < n\lambda$. Figure 4 shows the convergence of the speeds $v_i$, $i = 1,...,4$ to the desired speed $v^*$. Figure 5 shows the evolution of the intervehicle distances and their convergence to the spacing equilibrium $n^{-1}R = 32.5$ (recall the unique equilibrium point $E$ given by (4.15)). Finally, Figure 6 shows the acceleration $F_i(t)$ for each vehicle over time.

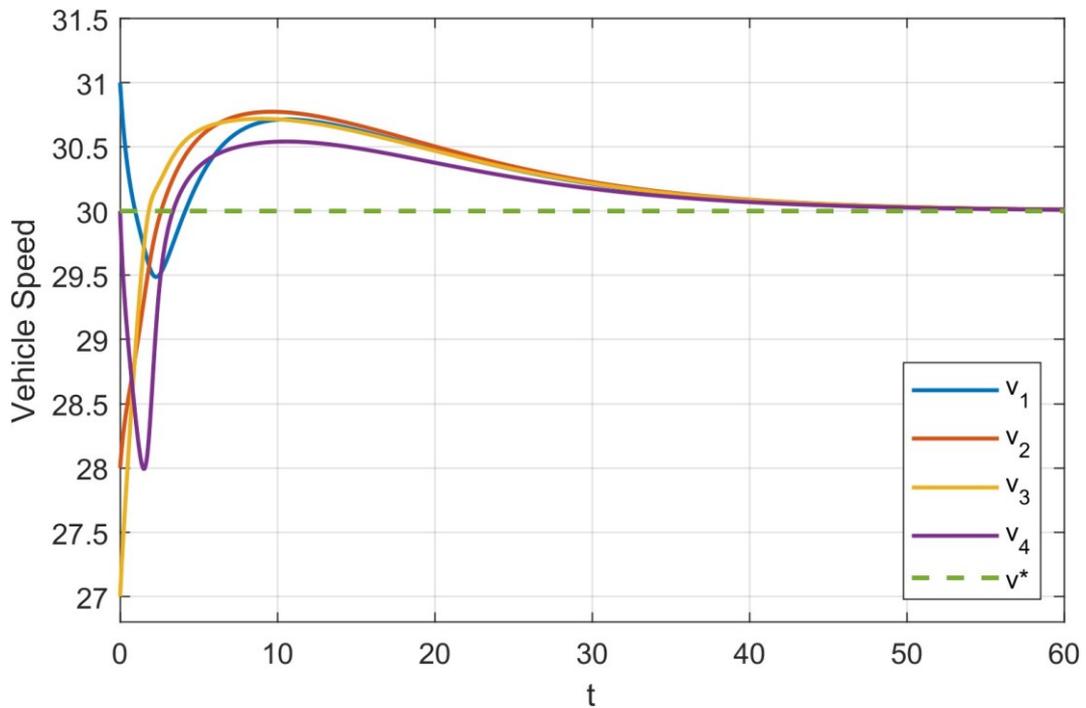

**Figure 4:** Convergence of vehicle speeds to the desired speed $v^*$ for the case $R < n\lambda$..

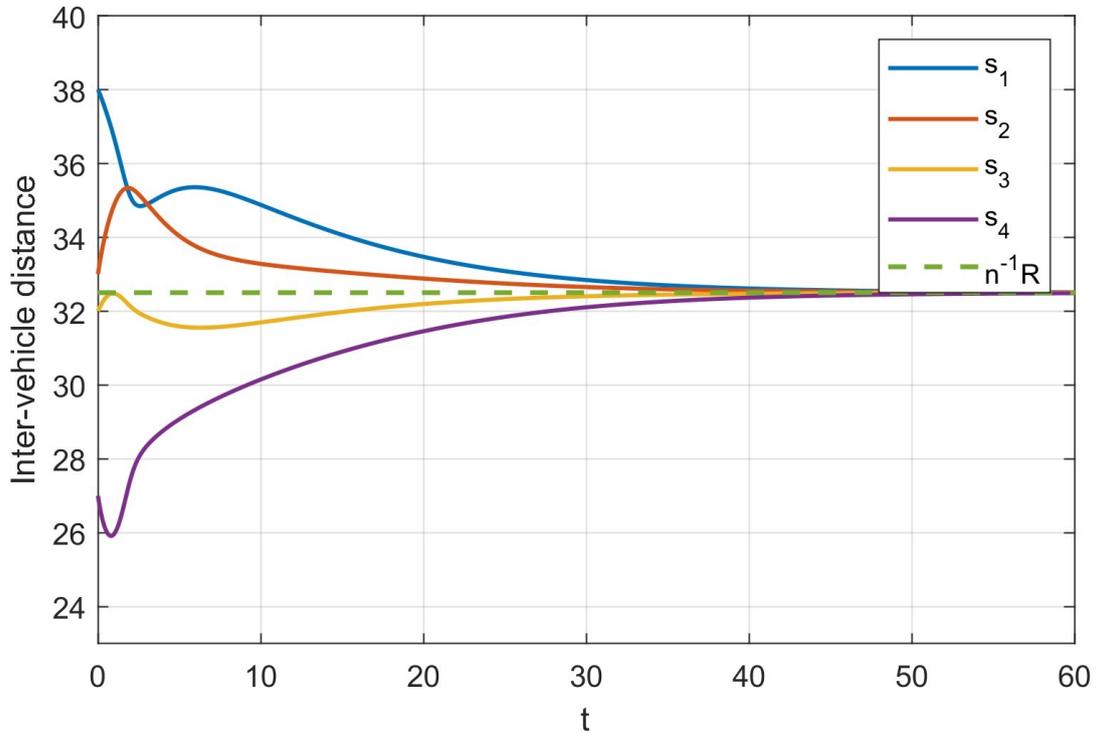

**Figure 5:** Evolution of inter-vehicle distance $s_i$ for the case $R < n\lambda$.

Finally, Figure 7, shows the convergence of the (normalized) logarithms of the Lyapunov functions $H(s,v)$ and $U(s,v)$ defined by (4.5) and (4.17), respectively. It can be seen that the convergence of $\log(U(s,v))$ is linear (or equivalently the convergence of $U(s,v)$ is exponential) as expected by Theorem 3, with convergence rate given by $c = \min\left(\mu, \vartheta^2 b'(0) V''\left(n^{-1}R\right)\mu_n\right)$, $\mu_n = 2$, $\vartheta = 0.9$ (shown with yellow in Figure 7).

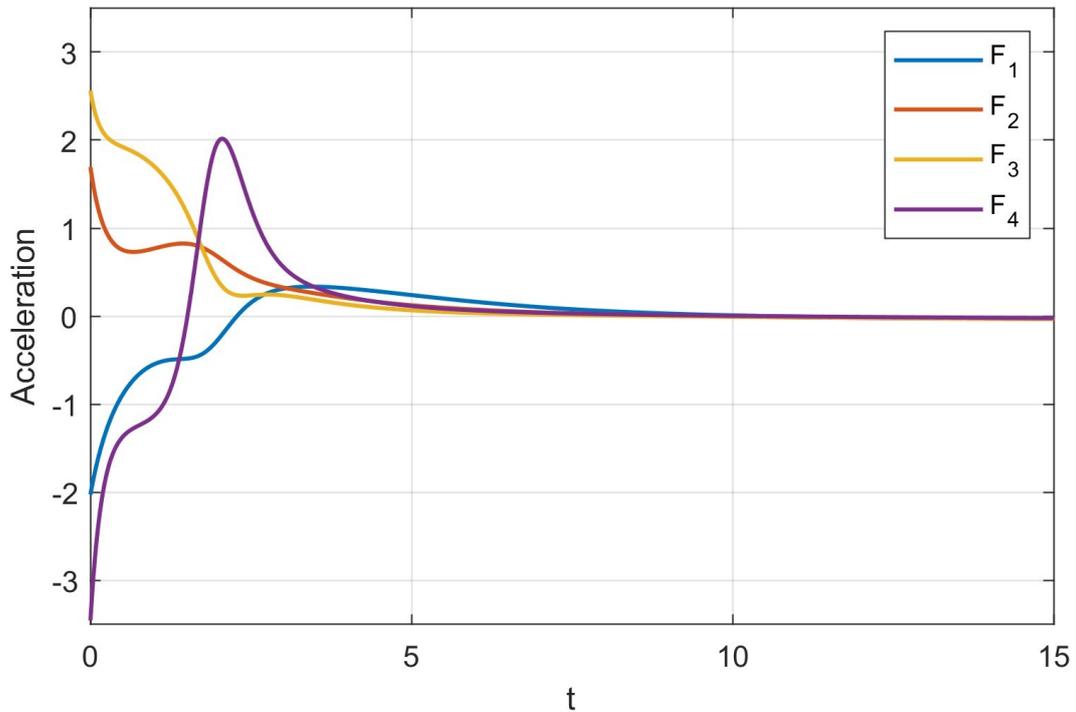

**Figure 6:** Convergence of vehicle accelerations $F_i$ for the case $R < n\lambda$.

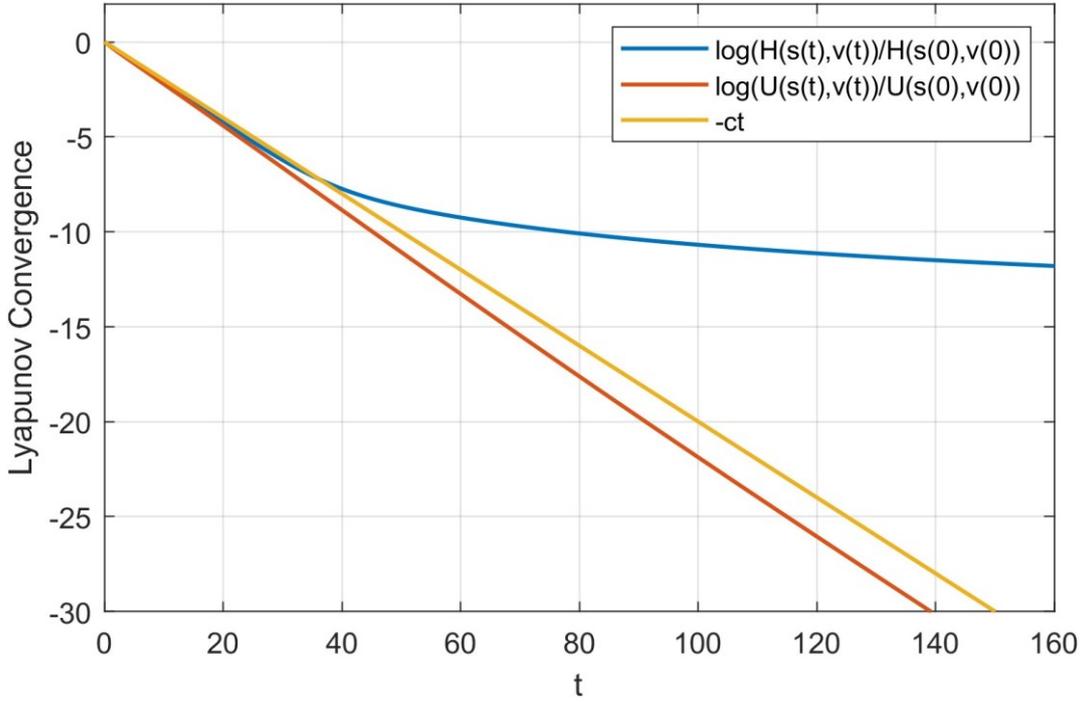

**Figure 7:** Convergence of the Lyapunov functions $H(s,v)$ and $U(s,v)$.

*5.2 Numerical Investigation of String Stability*

Here, we study numerically the string stability properties of the cruise controller (4.8) and its sensitivity to external disturbances on a ring-road with $R=130$. We consider $n=6$ vehicles with the speed of the first vehicle acting as a disturbance; namely, we consider the system

$$\dot{s}_i = v_{i-1} - v_i, \quad i=1,...,n$$
$$\dot{v}_i = F_i, \quad i=2,...,n$$

with $v_0 = v_n$ and

$$v_1(t) = \begin{cases} v^*, & t \in [0, \pi/2) \cup [5\pi/2, +\infty) \\ v^* + d\cos(t) & t \in [\pi/2, 5\pi/2) \end{cases},$$

where $d \in (0, \min(v^*, v_{max} - v^*))$ is constant such that $v_1(t) \in (0, v_{max})$ for all $t \geq 0$. We consider the following setting: $v_{max} = 35$, $v^* = 20$, $L=5$, $\mu = 0.1$, $\lambda = 40$, $d = 14$ with $V(s)$ and $b(s)$ given by (7.1) and (7.2). We further assume that all vehicles start at an equilibrium position $(n^{-1}R1_n, v^*1_n)$. Figure 8 shows the evolution of the speeds $v_i(t)$ of all vehicles when $v_1(t)$ acts as a disturbance. Moreover Figure 9 shows the speed deviation from the desired speed $|v_i(t) - v^*|$, $i=1,...,6$. It can be seen that during the time-interval $[\pi/2, \pi]$, when vehicle $i=1$ decelerates, the disturbance on the speed propagates backwards but weakens along the string of vehicles. On the other hand, when vehicle $i=1$ accelerates in the time interval $[\pi, 2\pi]$, and consequently affects more the vehicle in front ($i=6$) due to the bidirectional nature of the cruise controller and the geometry of the ring-road, the disturbances propagates forward to vehicle $i=6$, $i=5$, etc., but diminishes along the way.

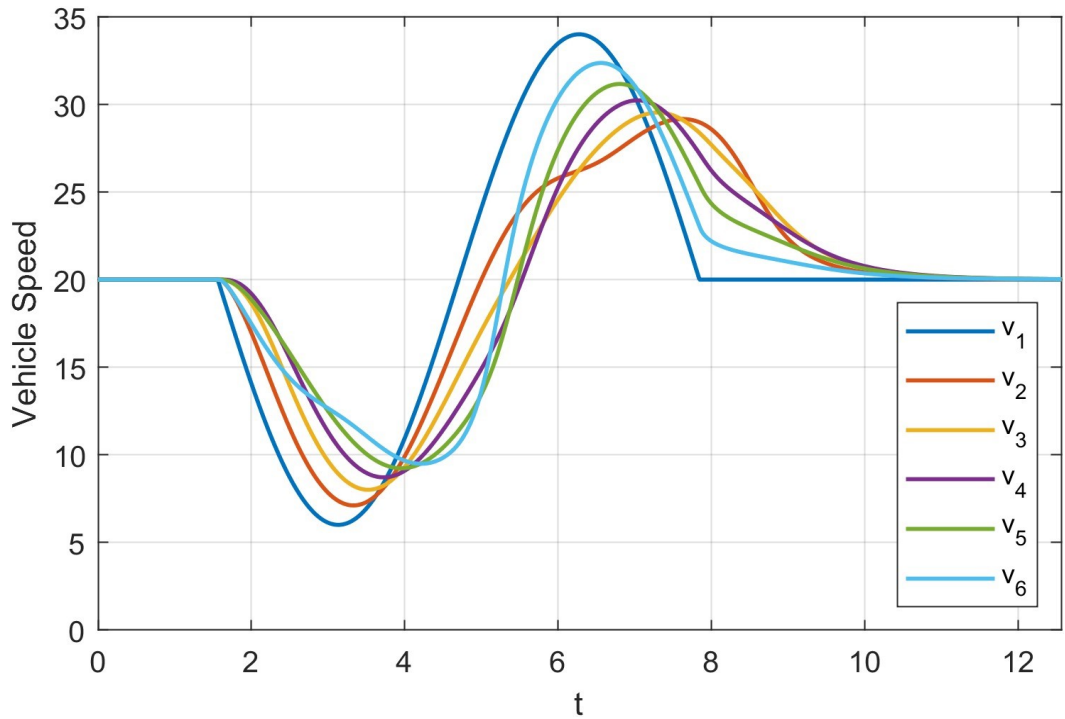

**Figure 8:** Evolution of speeds $v_i(t)$, $i = 1,...,6$ when $v_1(t)$ acts a disturbance.

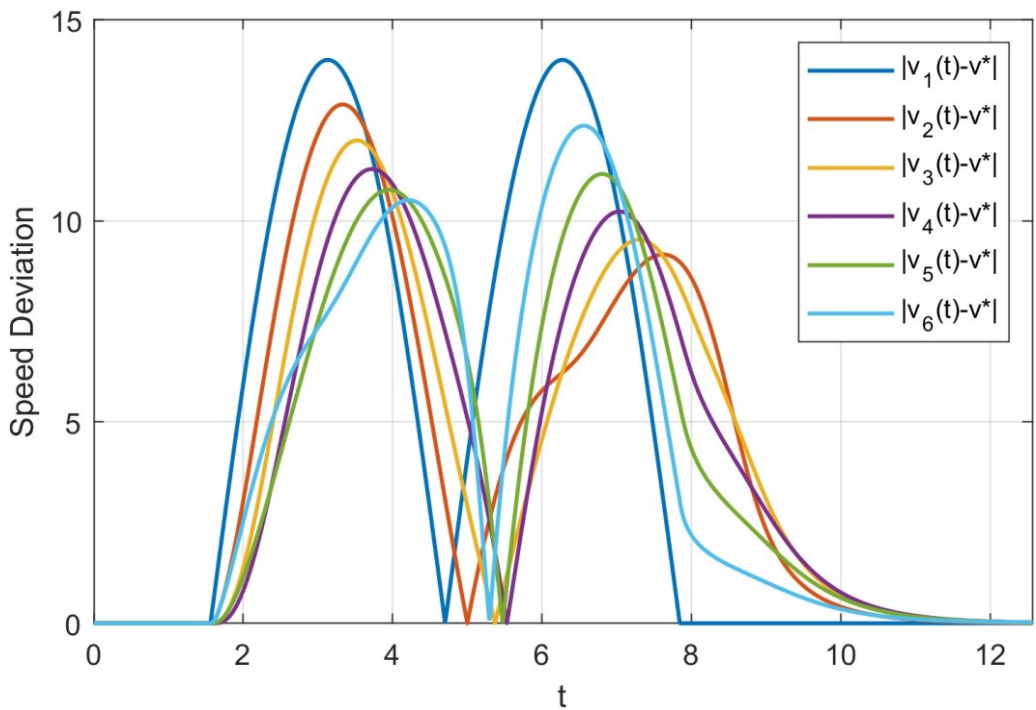

**Figure 9:** Evolution of the deviation of speeds $v_i(t)$ from the desired speed $v^*$.

## 5.3 Application on an Open Road.

In this section we demonstrate some properties of the cruise controller (4.8) and compare it with the cruise controller presented in [16]:

$$F_i = -k_i(s)(v_i - v^*) + V'(s_i) - V'(s_{i+1}) \quad (7.3)$$

where

$$k_i(s) = \tilde{\mu} + \frac{v_{max} g(V'(s_i) - V'(s_{i+1}))}{v^*(v_{max} - v^*)} - \frac{V'(s_i) - V'(s_{i+1})}{v^*}$$

$$g(x) = \frac{1}{2\varepsilon} \begin{cases} 0 & x \leq -\varepsilon \\ (x+\varepsilon)^2 & -\varepsilon < x < 0 \\ \varepsilon^2 + 2\varepsilon x & x \geq 0 \end{cases}$$

with $\varepsilon$ being a positive constant. It should be noted here that there are important differences between the controller (7.3) and the controller (4.8). The controller in (7.3) guarantees that $v_i(t) \in (0, v_{max})$, $i = 1,...,n$, for $t \geq 0$ due to the state-dependent controller gain $k_i(s)$, while the same property is guaranteed in controller (4.8) directly from the Lyapunov function $H_S$ in (5.2) which acts as a barrier function preventing the state to escape from the set $\Omega_S$. Another important difference is that in (7.3), the desired speed $v^*$ in the friction term $-k_i(s)(v_i - v^*)$ is constant, while in (4.8) it may decrease when the distance to the preceding and following vehicles is small; recall definition (4.6).

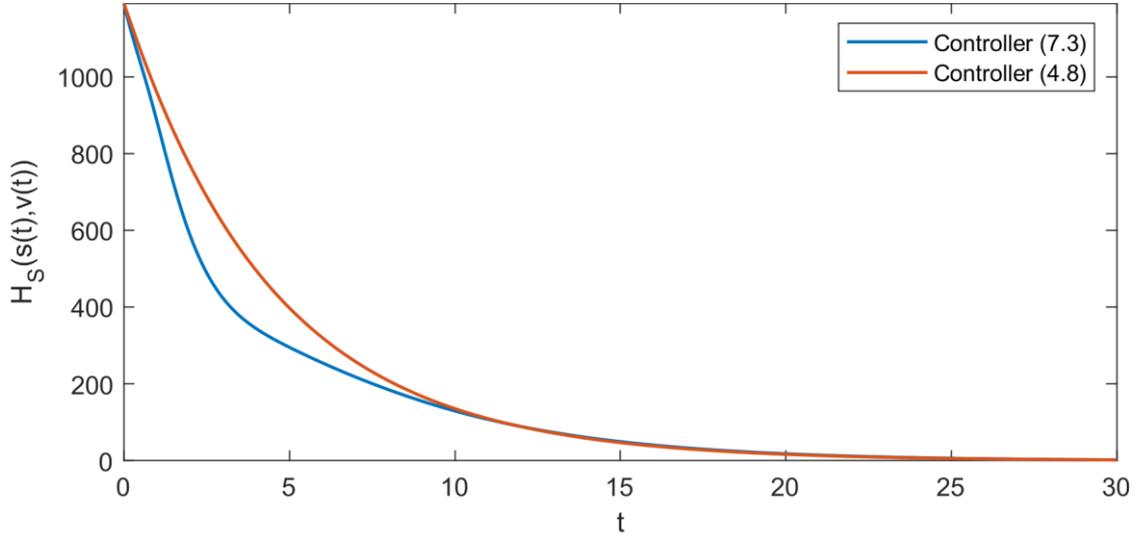

**Figure 10:** Convergence of the Lyapunov function $H_S(s(t), v(t))$ in (4.5) along the solution of the closed-loop systems (2.2), (4.8) and (2.2), (7.3).

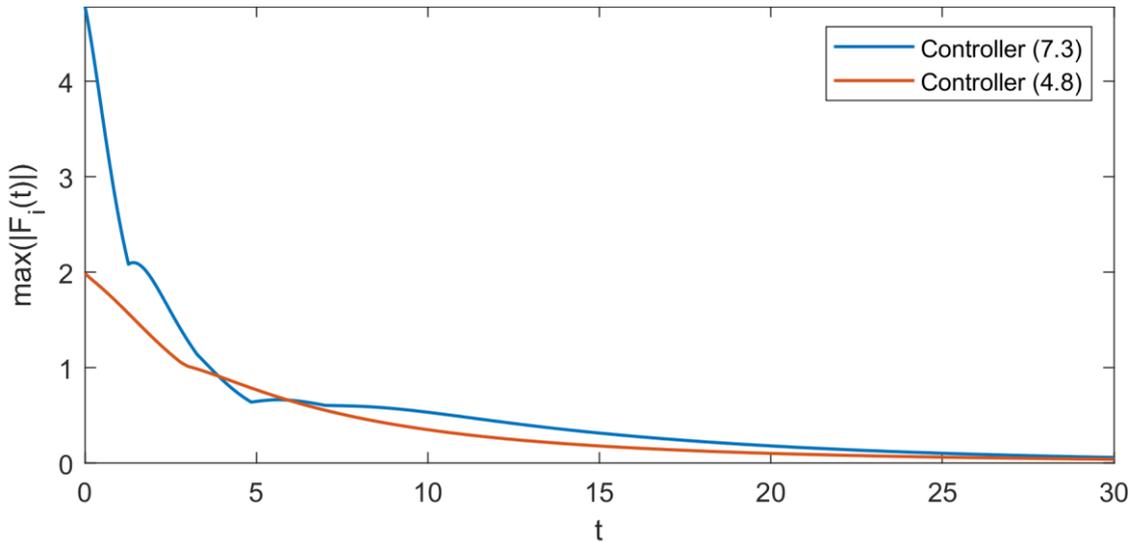

**Figure 11:** Evolution of $\max_{i=1,\ldots,n}\{|F_i(t)|\}$ for both controllers (4.8) and (7.3).

We consider a string of $n=5$ vehicles and use the following parameters and initial conditions for both controllers (4.8) and (7.3): $\lambda = 35$, $L = 5$, $v^* = 30$, $v_{\max} = 35$, $\varepsilon = 0.1$, $\mu = \tilde{\mu} = 0.1$, $q = (\lambda)^{-3}$ and $s_i(0) = 19$, $i = 2,\ldots,5$, $v_i(0) = 20$, $i = 1,\ldots,5$. To compare the controllers (4.8) and (7.3), we use the Lyapunov function $H_S(s,v)$ in (5.2) as a measure to examine their convergence to the set of equilibrium points $E_S$. Figure 10 shows that, with the above parameter choices, the convergence rate of the Lyapunov function to zero along the solution of each of the closed-loop systems (2.2), (4.8) and (2.2), (7.3) is similar for both controllers. On the other hand, Figure 11, which shows the maximum the acceleration of all vehicles over time, $\max_{i=1,\ldots,n}\{|F_i(t)|\}$, $t \geq 0$, demonstrates that the acceleration of controller (4.8) is significantly lower than the acceleration of the controller (7.3). This behaviour is due to the term $-\mu(v_i - f_i(s))$ in (4.8) which acts as a friction term that drives each vehicle to the spacing-dependent desired speed $f_i(s)$; while the friction term $k_i(s)(v_i - v^*)$ in (7.3) includes the state-dependent controller gain $k_i(s)$ whose value may increase for small inter-vehicle distances.

## 6. Conclusions

In this paper we presented a cruise controller with collision avoidance applied to the cases of a ring-road and an open road. The proposed controller is decentralized and uses only spacing and speed information from the preceding and following vehicles. For the case of the ring-road we distinguished two cases based on the interaction distance, the length of the ring-road as well as the number of vehicles, which determine whether we have a single equilibrium or a continuum of equilibrium points. In both cases we have shown $KL$ estimates that establish uniform convergence to the set of equilibria. In addition, when there exists a single equilibrium point, then every solution of the closed-loop system converges exponentially fast. Finally, we have studied the case of an open road and we have also provided certain $KL$ estimates for uniform convergence.

# Acknowledgments

The research leading to these results has received funding from the European Research Council under the European Union's Horizon 2020 Research and Innovation programme/ ERC Grant Agreement n. [833915], project TrafficFluid.